\documentclass[proc]{edpsmath}
\usepackage{xcolor}
\usepackage{amsmath}
\usepackage{graphicx}
\usepackage{subcaption}
\usepackage{booktabs}
\usepackage{makecell}
\usepackage{soul}

\usepackage{tikz}
\usetikzlibrary{positioning}

\def\eg{{e.g.,\ }}

\def\sindy{{SINDy}}
\definecolor{review1}{RGB}{117,112,179}

\usepackage[colorlinks=true]{hyperref}
\usepackage[nameinlink,capitalize,noabbrev]{cleveref}

\makeatletter
\newcommand{\sinname}{sin}

\newcommand{\sin@name}{%
  \smash{\operator@font\sinname}\vphantom{s}%
}
\DeclareRobustCommand{\sin}{%
  \mathop{\vphantom{\operator@font\sinname}}\!%
  \qopname\relax o{\sin@name}%
}
\makeatother

\numberwithin{equation}{section}

\begin{document}

\title{Generalizing the \sindy\ approach with nested neural networks}
%
\author{Camilla Fiorini}
\address{Conservatoire National des Arts et Métiers, Laboratoire M2N, 75003, Paris}
\author{Clément Flint}
\address{Université de Strasbourg, CNRS, Inria, IRMA, F-67000 Strasbourg, France}
\author{Louis Fostier}
\address{PRC, INRAE, CNRS, Université de Tours, Nouzilly, France}
\secondaddress{Université Paris-Saclay, Inria, Inria Saclay-Île-de-France, Palaiseau, France}
\author{Emmanuel Franck}
\sameaddress{2}
\author{Reyhaneh Hashemi}
\address{Aix--Marseille University, Aix-en-Provence, France}
\author{Victor Michel-Dansac}
\sameaddress{2}
\author{Wassim Tenachi}
\address{Universit\'e de Strasbourg, CNRS, Observatoire astronomique de Strasbourg, UMR 7550, F-67000 Strasbourg, France}
%
%
\begin{abstract}
    Symbolic Regression (SR) is a widely studied field of research that aims to infer symbolic expressions from data.
    A popular approach for SR is the Sparse Identification of Nonlinear Dynamical Systems (\sindy) framework, which uses sparse regression to identify governing equations from data.
    This study introduces an enhanced method, Nested SINDy, that aims to increase the expressivity of the \sindy\ approach thanks to a nested structure.
    Indeed, traditional symbolic regression and system identification methods often fail with complex systems that cannot be easily described analytically.
    Nested SINDy builds on the \sindy\ framework by introducing additional layers before and after the core \sindy\ layer.
    This allows the method to identify symbolic representations for a wider range of systems, including those with compositions and products of functions.
    We demonstrate the ability of the Nested SINDy approach to accurately find symbolic expressions for simple systems, such as basic trigonometric functions, and sparse (false but accurate) analytical representations for more complex systems.
    Our results highlight Nested SINDy's potential as a tool for symbolic regression, surpassing the traditional \sindy\ approach in terms of expressivity.
    However, we also note the challenges in the optimization process for Nested SINDy and suggest future research directions, including the designing of a more robust methodology for the optimization process.
    This study proves that Nested SINDy can effectively discover symbolic representations of dynamical systems from data, offering new opportunities for understanding complex systems through data-driven methods.
\end{abstract}
\begin{resume}
    La régression symbolique est un domaine de recherche bien établi qui cherche à déduire des expressions symboliques directement à partir de données.
    L'une des approches les plus reconnues dans ce domaine est celle de l'Identification Parcimonieuse de Systèmes Dynamiques non Linéaires (Sparse Identification of Nonlinear Dynamical Systems, SINDy), qui recourt à la régression parcimonieuse pour extraire les expressions des équations sous-tendant les données observées.
    Dans cette étude, nous proposons une méthode avancée, dénommée ``Nested SINDy'', qui vise à améliorer l'expressivité de la méthode SINDy grâce à une structure multi-couches.
    En effet, les limites des méthodes traditionnelles de régression symbolique et d'identification de systèmes concernent principalement les sytèmes complexes, dont la description analytique n'est pas aisée.
    En s'appuyant sur la méthode SINDy, Nested SINDy introduit des couches supplémentaires avant et après la couche principale de SINDy.
    Ceci lui permet d'inférer des représentations symboliques pour une plus grande gamme de systèmes, incluant ceux décrits par des équations formées de compositions et de produits de fonctions.
    Nous démontrons la capacité de Nested SINDy à identifier avec précision des expressions symboliques pour des équations simples, tels que des fonctions trigonométriques basiques, et à fournir des représentations analytiques creuses (fausses, mais précises) pour des systèmes plus complexes.
    Nos résultats mettent en lumière le potentiel de Nested SINDy comme outil puissant pour la régression symbolique, surpassant les capacités expressives de la méthode SINDy traditionnelle.
    Toutefois, nous soulignons également les défis rencontrés dans le processus d'optimisation de Nested SINDy et proposons des pistes pour des recherches futures, notamment en vue de développer une méthodologie d'optimisation plus robuste.
    Cette étude établit que Nested SINDy est capable de découvrir efficacement des représentations symboliques des systèmes dynamiques à partir de données, ouvrant ainsi de nouvelles voies pour la compréhension des systèmes complexes via des méthodes basées sur les données.
\end{resume}
\maketitle
\section{Introduction}


Symbolic regression (SR) consists in the inference of a free-form symbolic analytical function $f: \mathbb{R}^{n_1} \longrightarrow \mathbb{R}^{n_2}$ that fits $\mathbf{y} = f(\mathbf{x})$ given data $(\mathbf{x}, \mathbf{y})$. It is distinct from regular numerical optimization procedures in that it consists in a search in the space of functional forms themselves by optimizing the arrangement of mathematical symbols (\eg $+$, $-$, $\times$, $/$, $\sin$, $\cos$, $\exp$, $\log$, \dots).

The rationale for employing SR can be broadly categorized into the following three core objectives.
\begin{enumerate}
    \item SR can be used to produce models in the form of compact analytical expressions that are interpretable and intelligible. This objective is particularly vital in natural sciences, such as physics \cite{SR_phy_demo_TowardAIPhysicist}, where the capacity to explain phenomena is equally valuable as predictive prowess. This is typically probed by assessing the capability of a system to recover the exact symbolic functional form from its associated data. However, one should note that many SR approaches excelling in this metric are often bested in fit accuracy when exact symbolic recovery is unsuccessful \cite{SRBench}. In other contexts where the compactness and inherent intelligibility of expressions may not be as critical, significantly longer but more robustly accurate expressions ($> 10^3$ mathematical symbols) are desirable as SR still offers key advantages in such scenarios.
    \item SR demonstrates the advantage of producing models that frequently exhibit superior generalization properties when compared to neural networks. \cite{EQL_SymbolsInNN, Kamienny_EndToEndSR, Kamienny_SR_for_RL, SR_Generalization_SmallDataset, SR_of_controller_DSRbased}
    \item Another noteworthy advantage is the ability to create models that demand significantly fewer computational resources than extensive numerical models like neural networks. This efficiency becomes especially relevant in \textit{multi-query} scenarios such as control loops \cite{Kamienny_SR_for_RL, SR_of_controller_DSRbased}, optimization or uncertainty quantification, where models must be executed frequently, and thus computational efficiency is crucial.
\end{enumerate}

SR has traditionally been approached through genetic programming, where a population of candidate mathematical expressions undergoes iterative refinement using operations inspired by natural evolution, such as natural selection, crossover, and mutation. This approach includes well-known tools like Eureqa software \cite{EureqaPaper2009, EureqaPaper2011_AFP}, as well as more recent developments \cite{SBP_GP, GP_GOMEA, gplearn, OPERON}.  Additionally, SR has been explored using a diverse array of probabilistic methods \cite{Exhaustive_SR, Grammar_prior_wMC_SR, BSR_Bayesian_SR, Divide_and_conquer, GSR}. For recent SR reviews, refer to \cite{SRBench, SR_interpretability_review, SR_review_Angelis}.

The rise of neural networks and auto-differentiation\footnote{Leveraging the capabilities of deep learning libraries to meticulously track gradients associated with a set of parameters in relation to a numerical process, regardless of its intricacy.} has spurred significant efforts to incorporate these techniques into SR, challenging the dominance of Eureqa-like approaches \cite{SRBench, Benchmark_SR_phy, ExpressionSampler_benchmark}. Numerous methods for integrating neural networks into SR have been developed, ranging from advanced problem simplification schemes \cite{AIFeynman, AIFeynman2} to end-to-end supervised symbolic regression approaches in which neural networks are trained in a supervised manner to map datasets to their corresponding symbolic functions \cite{Kamienny_EndToEndSR, Biggio_EndToEndSR, DGSR, NeSymReS_EndToEndSR, SymFormer_EndToEndSR, Lalande_supervised, Boolformer, SNIP_supervised, Li2024GPT_supervised, MMSR_SR_multimodal_supervised, Botfip_supervised}. Unsupervised approaches also exist, where recurrent neural networks are trained through trial-and-error using reinforcement learning to generate analytical expressions that fit a given dataset \cite{PetersenDSR, SR_PG_improvements, uDSR_DSRbased, GFN_SR, Zheng_SymbolsInNN_DSRbased, ParetoFront_SR_RL, DSR_wTransformers_RL, SymQ_RL_SR, AlphaZero_SR_RL}.
Furthermore, it should be noted that it has been a major focus of the SR community to facilitate the incorporation of prior knowledge to constrain the search for functional forms by leveraging domain-specific knowledge \cite{PhySO_ApJ, PhySO_neurips, ClassSR, Bartlett_prior, DSR_priors, bayesian_prior, DSR_wikipedia_prior, DA_grammar_SR, PhySO_usage} and that \sindy-like frameworks as the one proposed here can accommodate such prior knowledge, as demonstrated in works like \cite{SISSO_DA}.

Supervised approaches offer rapid inference but lack a self-correction mechanism. If the generated expression is suboptimal, there are little means of correction. In contrast, unsupervised approaches enable iterative correction based on fit quality. However, they often rely on reinforcement learning frameworks to approximate gradients because direct optimization using auto-differentiation is infeasible due to the discrete nature of the problem, which involves discrete symbolic choices.

However, other unsupervised methods include neuro-symbolic approaches, wherein mathematical symbols are integrated into neural network frameworks. The goal being to sparsely fit the neural network to enable interpretability, generalization or even recover a compact mathematical expression. Prominent examples include \sindy~\cite{BruProKut2016}, which is central to this study, and others such as \cite{Martius_SymbolsInNN, ParFam, EQL_SymbolsInNN, Valle_SymbolsInNN, Kim_SymbolsInNN, Panju_SymbolsInNN, SISSO}.

\sindy-like approaches are the only type of unsupervised techniques capable of directly utilizing gradients from data to iteratively refine function expressions as they effectively render the discrete symbolic optimization problem continuous. Moreover, \sindy-like frameworks possess the advantage of being well-suited for exact symbolic recovery by enabling the creation of concise, intelligible analytical expressions through the promotion of sparse symbolic representations while yielding highly accurate and general expressions when exact symbolic recovery is unsuccessful or impossible. However, a limitation of the current \sindy\ framework is its inability to handle nested symbolic functions, which often results in suboptimal performances, especially in more complex problems as evidenced by comparative benchmarks (see \eg \cite{PhySO_ApJ}). This is the primary motivation for our study, where we introduce a Nested \sindy\ approach.

The paper is organized as follows.
First, we introduce the traditional \sindy\ approach in \cref{sec:sindy}.
Then, we present our Nested \sindy\ approach in \cref{sec:nested_sindy},
with two distinct architectures:
the PR model in \cref{sec:PR_model}
and the PRP model in \cref{sec:PRP_model},
the latter being more expressive, but also more challenging to train,
than the former.
The training procedure is explained in detail in \cref{sec:training},
and two main applications are tackled in \cref{sec:applications}:
function identification in \cref{subsec:function_discovery}
and ODE discovery in \cref{subsec:ode_discovery}.
Finally, we conclude in \cref{sec:conclusion}.

\section{The SINDy paradigm}
\label{sec:sindy}

In this section, we present the traditional \sindy\ approach,
as introduced in~\cite{BruProKut2016}.
The principle of the method is detailed in \cref{subsec:sindy_naive},
while the mathematical framework and notation are set in \cref{subsec:sindy_math}.

\subsection{Principle of the SINDy method}
\label{subsec:sindy_naive}

The Sparse Identification of Nonlinear Dynamical Systems
(\sindy) approach extends previous work in SR,
introducing innovations in sparse regression.
The \sindy\ approach seeks to deduce the governing equations
of a nonlinear dynamical system directly from observational data,
doing so in a concise and sparse way.
It is based on the essential assumption that
these governing equations can be succinctly expressed by
only a few significant terms, resulting in
a sparse representation within the space of
potential functions~\cite{BruProKut2016}.

More specifically, the \sindy\ approach involves
approximating a target function through a linear combination of
(potentially nonlinear) basis functions,
contained in a so-called library or dictionary $\mathcal{F}$.
For instance, $\mathcal{F}$ might include
constant, polynomial, or trigonometric functions:
\begin{equation}
    \label{eq:sindy_dico}
    \mathcal{F} =
    \begin{Bmatrix}
        x \mapsto 1,           &
        x \mapsto x,           &
        x \mapsto x^2,         &
        x \mapsto x^4,         &
        \dots                    \\
        x \mapsto \sin(x),     &
        x \mapsto \cos(x),     &
        x \mapsto \sin(2x),    &
        x \mapsto \cos(\pi x), &
        \dots
    \end{Bmatrix}.
\end{equation}
In order to achieve expressiveness, it is necessary to chose a large number of basis functions.
However, with that many basis functions,
there is a risk of losing interpretability.
To address this issue, a sparsity constraint is imposed
on the coefficients of the linear combination.

The main advantages of the \sindy\ method are the following:
\begin{itemize}
    \item The use of underlying convex optimization
          algorithms ensures the method's applicability to
          large-scale problems~\cite{BruProKut2016}.
    \item The resulting nonlinear model identification inherently
          balances model complexity
          (i.e., sparsity of the function to be learned)
          with accuracy, leading to strong generalization ability.
    \item \sindy\ automatically identifies the relevant terms
          in the dynamical system without making prior assumptions
          about the system's form, through the use of gradient descent.
\end{itemize}
The main limitations of the method include:
\begin{itemize}
    \item The necessity to carefully select the appropriate
          library $\mathcal{F}$ based on the available data:
          for instance, compositions or multiplications
          of simpler functions have to be included in the dictionary
          to correctly represent more complex target functions.
    \item Training (to determine which functions to keep
          in the dictionary) is more sensitive to
          initialization than with other approaches.
\end{itemize}

The goal of the following section is to set the mathematical
framework and the main notation associated to the \sindy\ method,
to be used throughout this paper.

\subsection{Mathematical framework and notation}
\label{subsec:sindy_math}

For the sake of simplicity,
we present the method in the case where the target function is
from $\mathbb{R}$ to $\mathbb{R}$,
but the approach can be extended to functions
from $\mathbb{R}^n$ to $\mathbb{R}^m$.

Given the data $(x_i, y_i)_{i=1,\ldots ,N }$,
we aim to find a function $f$ such that
$f(x_i) \approx y_i$ (which is nothing but a regression problem)
using the SINDy approach.

Let $\mathcal{F} = \left\{f_1, \dots, f_l \right\}$
be the aforementioned dictionary of basis functions.
For instance, it could be the one given by \eqref{eq:sindy_dico}.
We denote by $L(\mathcal{F}) = \text{Span}(\mathcal{F})$
the set of linear combinations of these basis functions,
defined by
\begin{equation}
    \label{eq:def_L_F}
    f \in L(\mathcal{F})
    \iff
    \exists \theta \in \mathbb{R}^l \text{ such that }
    f = \sum_{i=1}^l \theta_i f_i.
\end{equation}
The regression problem can then be formulated
as the following least squares problem:
\begin{equation*}
    \min_{ f \in L(\mathcal{F}) }
    \| Y - f(X) \|_2^2,
\end{equation*}
where $X = (x_1, \dots, x_N)^T$ and $Y = (y_1, \dots, y_N)^T$.
This problem can itself be reformulated in matrix form,
using the definition of the vector space $L(\mathcal{F})$:
\begin{equation}
    \label{eq:loss_sindy_no_lasso}
    \min_{ \theta \in \mathbb{R}^l }
    \| Y - \mathbb{F}(X) \theta \|_2^2
\end{equation}
where $
    [\mathbb{F}(X)]_{i,j} = (f_j(x_i)) \in
    \mathcal{M}_{N,l}(\mathbb{R})
$.

Numerous algorithms exist to solve this problem while
promoting sparsity. Without being exhaustive,
notable methods include the
standard STLSQ (sequentially thresholded least squares) and
the LARS (least-angle regression) methods.
Another approach involves adding a regularization term
on the coefficients of the linear combination to favor sparsity.
The most popular is Lasso regularization,
but others exist (SR3, SCAD, MCP,~\dots).
In this work, we focus on the Lasso approach.
Introducing a Lasso regularization term to promote sparsity,
the optimization problem for the SINDy approach becomes,
instead of~\eqref{eq:loss_sindy_no_lasso}:
\begin{equation}
    \label{eq:loss_sindy_lasso}
    \min_{ \theta \in \mathbb{R}^l }
    \| Y - \mathbb{F}(X) \theta \|_2^2 + \lambda \| \theta \|_1,
\end{equation}
where $\lambda > 0$ is a hyperparameter,
to be manually set when using the method.
The values of $\lambda$ will be reported when
using the method in the following sections.
Specialized optimization algorithms,
such as ADMM (alternating direction method of multipliers),
are effective in solving regression problems with regularization.

\section{The Nested SINDy approach}
\label{sec:nested_sindy}

As mentioned in \cref{subsec:sindy_naive},
one of the main limitations of the \sindy\ approach is
the choice of the nonlinear basis functions
populating the dictionary $\mathcal{F}$.
For instance, if the unknown function happens to be
a composition or a multiplication of simple functions,
we cannot find the correct expression unless
this specific composition or multiplication is in $\mathcal{F}$.
In this paper, we aim at relaxing this constraint by
introducing a way of composing simple functions,
without having to manually add
these compositions to the dictionary.

In the same spirit as the approach investigated
in~\cite{Martius_SymbolsInNN, EQL_SymbolsInNN},
we will enlarge the set $L(\mathcal{F})$.
We will proceed by analogy with a
standard approach in machine learning,
which involves considering models
with multi-layer neural networks rather than a single broad layer.
Here, instead of considering a single layer of
nonlinear functions, we explore an augmented architecture,
consisting of several such layers.
We will refer to this approach as Nested \sindy.

The cost to bear is the heightened complexity
of the optimization landscape.
Indeed, the optimization problem,
used to be the linear least squares problem
given by~\eqref{eq:loss_sindy_lasso}.
Now, it becomes a nonlinear problem,
since the matrix $\mathbb{F}(X)$
is replaced with a composition of nonlinear functions.
The new optimization problem
(still with Lasso regularization)
is formulated as follows:
\begin{equation*}
    \min_{\theta} \frac{1}{2} \lVert y - \mathcal{N}(x, \theta) \rVert_2^2 + \lambda \lVert \theta \rVert_1,
\end{equation*}
where $\mathcal{N}$ denotes our nested model,
parameterized by $\theta$.
Consequently, the resolution of this nonlinear optimization problem may be significantly more challenging and algorithms that work well in the linear case do not necessarily have the same performance in the nonlinear case.

The method then primarily depends on the choice of the
architecture $\mathcal{N}$.
The goal is to introduce new layers that achieve
favorable trade-offs between
expressivity and optimization complexity.
In this work, we propose two architectures,
which are described in the following sections:
the PR model (in \cref{sec:PR_model})
and the PRP model (in \cref{sec:PRP_model}).
Both of are based on the polynomial layer,
described in \cref{sec:poly_layer}.
From now on, we call the basic \sindy\ layer,
given by a projection~\eqref{eq:def_L_F} onto $L(\mathcal{F})$,
the radial layer.






\subsection{Polynomial layers}
\label{sec:poly_layer}

To improve the expressivity of the usual \sindy\ approach,
we introduce an additional layer,
to be applied after the usual \sindy\ expression.
This layer constructs a variety of monomials
from the input variable \( x \), represented as follows:
\[
    f_{\text{poly}}(x) = \sum_{i=0}^{d} \omega_i x^i,
\]
where \( d \) represents the maximum allowed polynomial degree
and \( \omega_i \) are the weights.
Note that \( \omega_0 \) is the constant part of the layer,
which corresponds to the bias in traditional neural networks.
For inputs with multiple variables,
the layer extends to a multivariate polynomial,
facilitating complex combinations of the variables.
For example, with two variables $x$ and $y$, we obtain
\[
    f_{\text{poly}}(x, y) = \sum_{i=0}^{d} \sum_{j=0}^{d} \omega_{i,j} x^i y^j,
\]
with \( \omega_{0,0} \) acting as the constant term
of the polynomial, effectively substituting the bias.
Another choice is to limit the sum over $i+j$
to a maximum value $d$ to reduce
the number of terms in the polynomial,
and thus obtain bivariate polynomials up to degree $d$.
This strategy is applied in \cref{sec:PRP_model};
with two variables $x$ and $y$, such a polynomial layer reads:
\begin{equation}
    \label{eq:poly_layer_fixed_degree}
    f_{\text{poly}}(x, y) = \sum_{0 \le i+j \le d }  \omega_{i,j} x^i y^j.
\end{equation}

\subsection{The PR Model}
\label{sec:PR_model}

The PR (Polynomial-Radial) model augments the
basic SINDy framework by introducing a polynomial layer
that operates before the usual \sindy\ radial layer.
As an example, for one input variable,
if $\mathcal{F} = \{ \sin, \cos \}$ and $d = 2$,
then the PR model can learn all functions with expression
\begin{equation*}
    \lambda \cos(a_1 + b_1 x + c_1 x^2) +
    \mu \sin(a_2 + b_2 x + c_2 x^2).
\end{equation*}
This is way more expressive than standard \sindy,
where functions such as $x \mapsto \cos(2x)$
or $x \mapsto \sin(1 + x^2)$
would have to be manually added to the dictionary.

\begin{figure}[!ht]
    \centering
    \begin{tikzpicture}
        \tikzset{
            box/.style={draw, thick, rectangle},
            cell/.style={draw, rectangle, minimum height=1cm, minimum width=1cm}
        }

        \node[cell] (x_in) at (0,-1cm) {$x$};
        \node[cell] (y_in) at (0,-2cm) {$y$};
        \node[inner sep=0pt, outer sep=0pt, right=0cm of x_in] (input) {};

        \node[cell] (x_poly) at (3cm,0.5cm) {$x$};
        \node[cell] (xy_poly) at (3cm,-0.5cm) {$xy$};
        \node[cell] (xy2_poly) at (3cm,-1.5cm) {$xy^2$};
        \node[cell] (dots_poly) at (3cm,-2.5cm) {$\vdots$};
        \node[cell] (y2_poly) at (3cm,-3.5cm) {$y^2$};
        \node[inner sep=0pt, outer sep=0pt, right=0cm of y2_poly] (polynomial) {};

        \node[cell, minimum height=4cm] (linear_cell) at (6cm,-1.5cm) {};
        \node[inner sep=0pt, outer sep=0pt, right=0cm of linear_cell] (linear) {};

        \node[cell, minimum height=4cm] (radial) at (9cm,-1.5cm) {};
        \node[inner sep=0pt, outer sep=0pt, right=0cm of radial] (radial) {};

        \node[cell] (output) at (12cm,-1.5cm) {output};

        \draw[->] (x_in) -- (x_poly);
        \draw[->] (x_in) -- (xy_poly);
        \draw[->] (x_in) -- (xy2_poly);
        \draw[->] (y_in) -- (xy_poly);
        \draw[->] (y_in) -- (xy2_poly);
        \draw[->] (y_in) -- (y2_poly);

        \draw[->] (x_poly) -- (linear_cell);
        \draw[->] (xy_poly) -- (linear_cell);
        \draw[->] (xy2_poly) -- (linear_cell);
        \draw[->] (y2_poly) -- (linear_cell);
        \draw[->] (x_poly) -- (5.5cm,0cm);
        \draw[->] (x_poly) -- (5.5cm,-2cm);
        \draw[->] (xy_poly) -- (5.5cm,-0.5cm);
        \draw[->] (dots_poly) -- (5.5cm,-2.5cm);
        \draw[->] (y2_poly) -- (5.5cm,-3cm);
        \draw[->] (y2_poly) -- (5.5cm,-1cm);

        \draw[->] (6.5cm,0cm) -- (8.5cm,0cm);
        \draw[->] (6.5cm,-0.5cm) -- (8.5cm,-0.5cm);
        \draw[->] (6.5cm,-1cm) -- (8.5cm,-1cm);
        \draw[->] (6.5cm,-1.5cm) -- (8.5cm,-1.5cm);
        \draw[->] (6.5cm,-2cm) -- (8.5cm,-2cm);
        \draw[->] (6.5cm,-2.5cm) -- (8.5cm,-2.5cm);
        \draw[->] (6.5cm,-3cm) -- (8.5cm,-3cm);

        \draw[->] (9.5cm,0cm) -- (output);
        \draw[->] (9.5cm,-0.5cm) -- (output);
        \draw[->] (9.5cm,-1cm) -- (output);
        \draw[->] (9.5cm,-1.5cm) -- (output);
        \draw[->] (9.5cm,-2cm) -- (output);
        \draw[->] (9.5cm,-2.5cm) -- (output);
        \draw[->] (9.5cm,-3cm) -- (output);

        \node[inner sep=0pt] (input_label)at (0,-4.5cm) {Input};
        \node[inner sep=0pt] (polynomial_label) at (3cm,-4.5cm) {Polynomial};
        \node[inner sep=0pt] (linear_label) at (6cm,-4.5cm) {Linear};
        \node[inner sep=0pt] (radial_label) at (9cm,-4.5cm) {Radial};
        \node[inner sep=0pt] (output_label) at (12cm,-4.5cm) {Output};

    \end{tikzpicture}
    \caption{Structure of the PR model for $d=2$ and 2 input variables.}
    \label{fig:pr_model}
\end{figure}
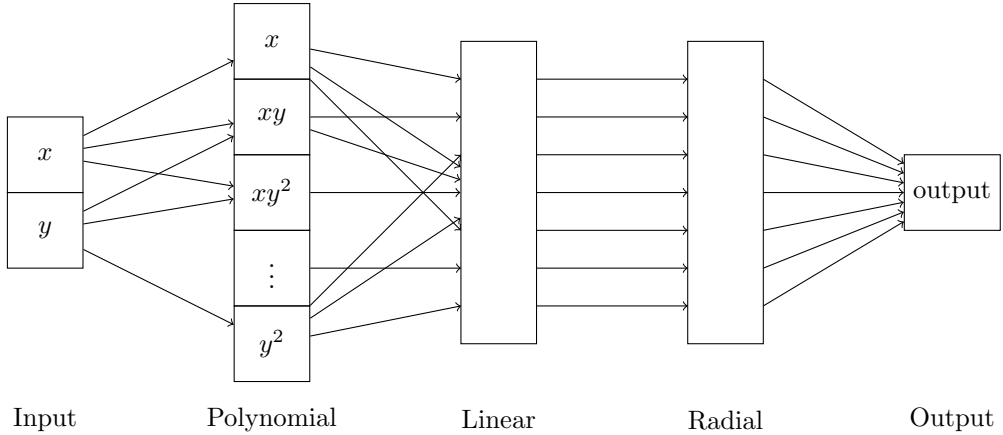

This PR layer can also be seen as a pure polynomial
layer combined with a linear layer.
We, hence, consider the PR model to have four layers:
a polynomial layer from \cref{sec:poly_layer},
a (fully-connected) linear layer, a radial layer,
and a final linear layer.
The last two layers are identical to the standard SINDy model,
while the first two layers are new additions.
\cref{fig:pr_model} illustrates this structure.

The full expression of the model is:
\begin{equation}
    \label{eq:pr_model}
    f_{\theta,\text{PR}}(x) = \sum\limits_{j=1}^l c_j f_j\left(\sum\limits_{i=0}^{d}\omega_{i,j}x^i\right) + \beta,
\end{equation}
where \( \theta \) includes
all the trainable parameters of the model.
The functions \( f_j \) are derived from
the specified function set \( \mathcal{F} \),
while \( \omega_{i,j} \) are the weights of the polynomial layer,
and \( c_j \) and \( \beta \) are the weights and bias
of the final linear layer, respectively.
The parameters \( \omega_{i,j} \) now depend on \( j \)
since they are in the \( j \)-th radial layer.
It is important to note that the radial and polynomial layers
are not associated with any adjustable parameters.

The main advantage of the PR model is its enhanced expressivity,
which facilitates the creation of linear combinations
both before and after the radial layer.
The polynomial layer allows the model to identify
more complex functions and to integrate various inputs
effectively in the case of multivariate data.
Compared with the traditional SINDy approach,
the model benefits from a reduced need for
an extensive dictionary because it is capable of
discovering linear combinations of the functions contained
within the dictionary.
Observations from our experiments suggest that
the training process of the model is capable of converging
to correct solutions even for non-trivial problems.
This will be highlighted in~\cref{sec:applications}.

\subsection{The PRP Model}
\label{sec:PRP_model}

The PRP (Polynomial-Radial-Polynomial) model enhances
the PR model by introducing an additional polynomial layer
following the radial layer.
This structure significantly improves the ability of the model
to represent complex interactions within datasets.
For example, the PRP model can express
the function $f(x) = \arctan(x) \sin(x)$,
assuming arctangent and sine are in the dictionary,
while the PR model cannot, as it is not a linear combination of the functions contained
within the dictionary.
In the PRP model, the outputs of the radial layer are
first processed through a fully-connected linear layer.
This linear layer serves as an intermediate stage,
transforming the outputs of the radial layer into
a new set of variables.
These variables are then fed into a subsequent polynomial layer,
which allows for the formation of various
monomial combinations of the outputs of the linear layers,
such as squaring or multiplying them.

\begin{figure}[htbp]
    \centering
    \begin{tikzpicture}
        \tikzset{
            box/.style={draw, thick, rectangle},
            cell/.style={draw, rectangle, minimum height=1cm, minimum width=1cm}
        }

        \node[cell] (x_in) at (0,-1cm) {$x$};
        \node[cell] (y_in) at (0,-2cm) {$y$};
        \node[inner sep=0pt, outer sep=0pt, right=0cm of x_in] (input) {};

        \node[cell] (x_poly) at (2.5cm,0.5cm) {$x$};
        \node[cell] (xy_poly) at (2.5cm,-0.5cm) {$xy$};
        \node[cell] (xy2_poly) at (2.5cm,-1.5cm) {$xy^2$};
        \node[cell] (dots_poly) at (2.5cm,-2.5cm) {$\vdots$};
        \node[cell] (y2_poly) at (2.5cm,-3.5cm) {$y^2$};
        \node[inner sep=0pt, outer sep=0pt, right=0cm of y2_poly] (polynomial) {};

        \node[cell, minimum height=4cm] (linear_cell) at (5cm,-1.5cm) {};
        \node[inner sep=0pt, outer sep=0pt, right=0cm of linear_cell] (linear) {};

        \node[cell, minimum height=4cm] (radial) at (7.5cm,-1.5cm) {};
        \node[inner sep=0pt, outer sep=0pt, right=0cm of radial] (radial) {};

        \node[cell] (linear_2_cell_1) at (10cm,-1cm) {$x_1$};
        \node[cell] (linear_2_cell_2) at (10cm,-2cm) {$x_2$};
        \node[inner sep=0pt, outer sep=0pt, right=0cm of linear_2_cell_1] (linear_2) {};

        \node[cell] (x1_poly_2) at (12.5cm,0.5cm) {$x_1$};
        \node[cell] (x1x2_poly_2) at (12.5cm,-0.5cm) {$x_1x_2$};
        \node[cell] (x1x22_poly_2) at (12.5cm,-1.5cm) {$x_1x_2^2$};
        \node[cell] (dots_poly_2) at (12.5cm,-2.5cm) {$\vdots$};
        \node[cell] (x2_poly_2) at (12.5cm,-3.5cm) {$x_2$};

        \node[cell] (output) at (15cm,-1.5cm) {output};

        \draw[->] (x_in) -- (x_poly);
        \draw[->] (x_in) -- (xy_poly);
        \draw[->] (x_in) -- (xy2_poly);
        \draw[->] (y_in) -- (xy_poly);
        \draw[->] (y_in) -- (xy2_poly);
        \draw[->] (y_in) -- (y2_poly);

        \draw[->] (x_poly) -- (linear_cell);
        \draw[->] (xy_poly) -- (linear_cell);
        \draw[->] (xy2_poly) -- (linear_cell);
        \draw[->] (y2_poly) -- (linear_cell);
        \draw[->] (x_poly) -- (4.5cm,0cm);
        \draw[->] (xy_poly) -- (4.5cm,-0.5cm);
        \draw[->] (dots_poly) -- (4.5cm,-2.5cm);
        \draw[->] (y2_poly) -- (4.5cm,-3cm);

        \draw[->] (5.5cm,0cm) -- (7cm,0cm);
        \draw[->] (5.5cm,-0.5cm) -- (7cm,-0.5cm);
        \draw[->] (5.5cm,-1cm) -- (7cm,-1cm);
        \draw[->] (5.5cm,-1.5cm) -- (7cm,-1.5cm);
        \draw[->] (5.5cm,-2cm) -- (7cm,-2cm);
        \draw[->] (5.5cm,-2.5cm) -- (7cm,-2.5cm);
        \draw[->] (5.5cm,-3cm) -- (7cm,-3cm);

        \draw[->] (8cm,0cm) -- (linear_2_cell_1);
        \draw[->] (8cm,-0.5cm) -- (linear_2_cell_1);
        \draw[->] (8cm,-1cm) -- (linear_2_cell_1);
        \draw[->] (8cm,-1.5cm) -- (linear_2_cell_1);
        \draw[->] (8cm,-2cm) -- (linear_2_cell_1);
        \draw[->] (8cm,-2.5cm) -- (linear_2_cell_1);
        \draw[->] (8cm,-3cm) -- (linear_2_cell_1);
        \draw[->] (8cm,0cm) -- (linear_2_cell_2);
        \draw[->] (8cm,-0.5cm) -- (linear_2_cell_2);
        \draw[->] (8cm,-1cm) -- (linear_2_cell_2);
        \draw[->] (8cm,-1.5cm) -- (linear_2_cell_2);
        \draw[->] (8cm,-2cm) -- (linear_2_cell_2);
        \draw[->] (8cm,-2.5cm) -- (linear_2_cell_2);
        \draw[->] (8cm,-3cm) -- (linear_2_cell_2);

        \draw[->] (linear_2_cell_1) -- (x1_poly_2);
        \draw[->] (linear_2_cell_1) -- (x1x2_poly_2);
        \draw[->] (linear_2_cell_1) -- (x1x22_poly_2);
        \draw[->] (linear_2_cell_2) -- (x1x2_poly_2);
        \draw[->] (linear_2_cell_2) -- (x1x22_poly_2);
        \draw[->] (linear_2_cell_2) -- (x2_poly_2);

        \draw[->] (13cm,0.5cm) -- (output);
        \draw[->] (13cm,-0.5cm) -- (output);
        \draw[->] (13cm,-1.5cm) -- (output);
        \draw[->] (13cm,-3.5cm) -- (output);

        \node[inner sep=0pt] (input_label)at (0,-4.5cm) {Input};
        \node[inner sep=0pt] (polynomial_label) at (2.5cm,-4.5cm) {Polynomial};
        \node[inner sep=0pt] (linear_label) at (5cm,-4.5cm) {Linear};
        \node[inner sep=0pt] (radial_label) at (7.5cm,-4.5cm) {Radial};
        \node[inner sep=0pt] (linear_2_label) at (10cm,-4.5cm) {Linear};
        \node[inner sep=0pt] (polynomial_2_label) at (12.5cm,-4.5cm) {Polynomial};
        \node[inner sep=0pt] (output_label) at (15cm,-4.5cm) {Output};

    \end{tikzpicture}
    \caption{Structure of the PRP model for $k=2$, $d=2$, and 2 input dimensions.}
    \label{fig:prp_model}
\end{figure}
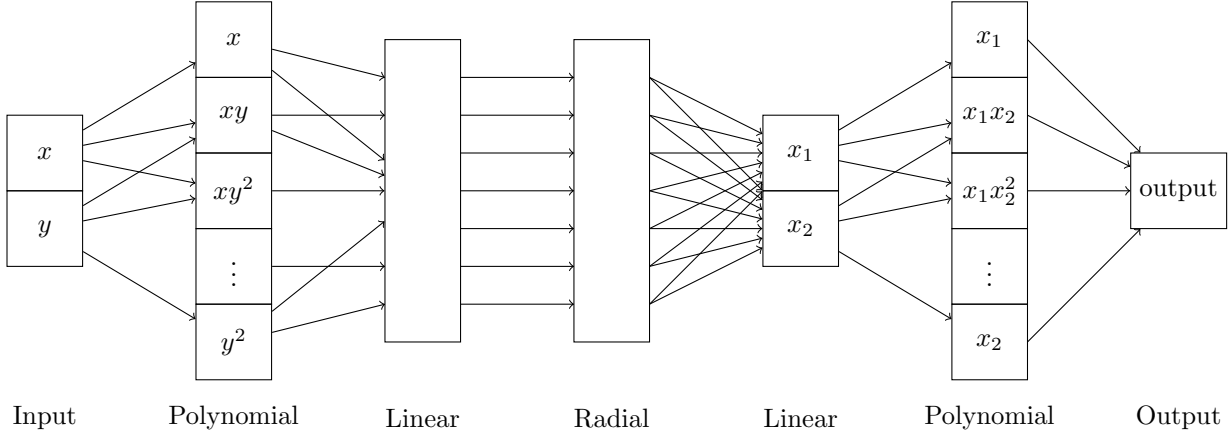

The mathematical expression of the PRP model is given as:
\begin{equation}
    \label{eq:prp_model}
    f_{\theta,\text{PRP}}(x) = \sum\limits_{1 \le \lvert i \rvert \le d } \omega_{i_1, i_2, \dots, i_l}^{\text{PR}}f_{1,\text{PR}}(x)^{i_1} f_{2,\text{PR}}(x)^{i_2} \dots f_{k,\text{PR}}(x)^{i_k} + \beta',
\end{equation}
where $\lvert i \rvert = i_1 + i_2 + \hdots + i_k $
is the length of the multi-index $(i_1,\hdots,i_k)$,
\( \theta \) includes all the trainable parameters of the model,
\( \omega_{i_1, i_2, \dots, i_k}^{\text{PR}} \) are
the weights of the final linear layer,
\( \beta' \) is the bias of the final linear layer,
$k$ is the size of the output chosen for
the intermediate linear layer (set to $k=2$ in our experiments),
and \( f_{1,\text{PR}}(x) \),
\( f_{2,\text{PR}}(x) \), \dots,
\( f_{k,\text{PR}}(x) \) are given by the following equation:
\begin{equation*}
    f_{{k'},\text{PR}}(x) = \sum\limits_{j=1}^l c_{j,k'} f_j\left(\sum\limits_{i=0}^{d}\omega_{i,j}x^i\right) + \beta_{k'}
    \text{, \quad for all } 1 \le k' \le k.
\end{equation*}
\cref{fig:prp_model} shows a graphical representation
of this model.

The addition of a second polynomial layer in the PRP model
markedly increases the expressivity of the model.
It enables the model to capture more intricate relationships
in the data, particularly beneficial for complex datasets
where simpler models may fall short.
Therefore, the PRP model is especially good at handling
datasets with intricate variable interactions.

In summary, the PRP model, with its dual polynomial layers,
offers a sophisticated extension of the SINDy approach.
It provides a powerful framework for modeling complex systems,
capable of capturing higher-order interactions and
nonlinear relationships inherent in the data.

\subsection{Downsides of the Nested SINDy approach}

The Nested SINDy approach complexifies the optimization
landscape, which destabilizes the training process.
We can reasonably assume that the addition of linear layers
creates local minima because of
their composition with the nonlinear layers.
Moreover, it adds a substantial amount of trainable parameters.
This is because the number of parameters
in the linear layers is quadratic,
whereas those in the \sindy\ approach grow linearly
with the number of functions in the dictionary.
This deviates from the original \sindy\ approach,
which is oriented towards function discovery,
and assumes that the dictionary can become arbitrarily large.
In the Nested SINDy approach,
the number of functions in the dictionary should remain
relatively small to avoid an explosion in the learning time.
The next section addresses some training-related issues.

\section{Training the nested \sindy\ model}
\label{sec:training}





We experimented and combined various strategies to best
overcome challenges associated with nonlinear optimization.
These strategies are given below,
where we mention how we actually used them for training.
Values of the hyperparameters introduced
in this section will be given in
\cref{sec:applications}.

\subsection{Adding a regularization term to enforce sparsity}

In our framework, the Lasso regularization term is added to the loss function to enforce sparsity in the model.
We have tested different strategies to adapt the Lasso coefficient during training:
\begin{itemize}
    \item A constant Lasso parameter throughout the training,
          which is the standard approach.
    \item A varying Lasso coefficient,
          initially set to zero for the early epochs,
          and then taken oscillating around some constant value,
          depending on the epoch.
          The intuition of this idea is that ``shaking'' the learning landscape
          helps to get out of local minima, in a direction that is still
          relevant to one of the two objectives (sparsity versus mean squared error).
    \item We also tried selecting neuron-dependent
          Lasso coefficients to promote
          specific functions or layers over others
          in the radial layer.
\end{itemize}
The approach that gave the most consistent results across the tested cases was to change the weight of the Lasso coefficient throughout the learning process. In the experiments, the Lasso coefficient is given by
\begin{equation}
    \lambda(\text{epoch}) = \lambda_0 (1 + 0.5 \sin(\text{epoch} / 10)),
    \label{eq : lasso oscillations}
\end{equation}
where $\lambda_0$ is the initial Lasso coefficient, and where the sine function is used to make the coefficient oscillate between $0.5 \lambda_0$ and $1.5 \lambda_0$.

\subsection{Pruning to enforce sparsity}
To enhance model sparsity, a complementary approach to Lasso regularization is to prune the neural network during training, reducing the number of parameters.
In order to prune the neural network, we remove a parameter $\theta_i$ from the set of parameters $(\theta_i)_i$ if the following two conditions are met:
\begin{itemize}
    \item The mean squared error (MSE) is below a predefined threshold value $\text{MSE}_\text{prune}$, and
    \item  $|\theta |$ is below a threshold value $\varepsilon_\text{prune}$ for a given number of epochs $n_\text{prune}$.
\end{itemize}

\subsection{Choosing the optimization algorithm}

We tested training our network with fairly standard optimization algorithms implemented in PyTorch but not specifically tailored to our problem: Adam, stochastic gradient descent, Limited-memory Broyden-Fletcher-Goldfarb-Shanno (LBFGS) algorithm.
For instance, LBFGS \cite{LBFGS} is uncommon in neural network training, as it was designed for optimizing constants in equations.
However, it works well in our case,
possibly because our optimization problem is close to a classical regression problem.
Moreover, it would be interesting to implement a more specific optimization algorithm that takes into account the form of our objective function (mean squared error plus regularization term), such as an ADMM algorithm coupled with a standard PyTorch optimizer. Another promising avenue for \sindy-like approaches is the \textit{basin-hopping} algorithm which combines LBFGS with global search techniques in order to avoid local minima as proposed in \cite{ParFam}.

\subsection{Adding noise to the gradient of the loss function}

During a training step, we can add random noise to the gradient just before updating weights.
This can be done at each training step, or only when the loss function does not vary sufficiently, e.g. when stuck in a local minimum.

The noise amplitude has to be well-tuned. It depends on the learning rate $\ell$ of the optimizer, the values of the parameters $\theta$, and on the actual mean squared error $L_\text{MSE}$:
\begin{equation}
    \nabla_\theta L \leftarrow \nabla_\theta L + \varepsilon (\ell,L_\text{MSE},\theta).
    \label{eq:add_noise_to_grad}
\end{equation}
with
\begin{equation*}
    \varepsilon \sim \mathcal{N} \left(0,\alpha \times \ell \times  \min (L_{\text{MSE}},1) \times \lvert \theta \lvert \right), \text{ and $\alpha$ a parameter to tune.}
\end{equation*}
By making the learning process less deterministic,
we hope to more easily escape local minima.
In the same spirit, we could also directly add noise to the
parameters when the loss function is stuck in a local minimum.

\subsection{Initializing the network parameters}

The training process is very sensitive to the weight initialization, given the presence of multiple local minima in the objective function. Furthermore, concerning the exploration of Ordinary Differential Equations (ODEs), for certain parameter values, the solution to the ODE may be ill-defined (or only locally defined), emphasizing the importance of careful weight initialization for our model. This is indeed one of the main limitations of the model and, consequently, an area for possible improvement.

\section{Applications}
\label{sec:applications}

In this work, we consider two main applications:
function discovery in \cref{subsec:function_discovery},
and Ordinary Differential Equation (ODE) discovery in \cref{subsec:ode_discovery}.
Function discovery consists in recovering the expression of a function from a dataset,
while ODE discovery aims at finding the differential equation
that governs the evolution of a system from a dataset.
Therefore, function discovery can be seen as a stepping stone towards ODE discovery.

\subsection{Function discovery}
\label{subsec:function_discovery}

We first tackle function discovery.
To that end, we present four test cases, of increasing complexity.

\subsubsection{Case 1: trigonometric function involving
    composition using the PR block}
\label{sec:function_case_1}

In \cref{subsec:sindy_naive}, we recalled that the standard \sindy\ method encounters difficulties when the target is a function defined as the composition of several simpler functions, unless that specific composite function is included in the dictionary. We wish to demonstrate, in such cases, the capability of the proposed PR nested \sindy\ method \eqref{eq:pr_model} described in \cref{sec:PR_model}.
To that end, the function $f: x \mapsto \cos(x^2)$ is considered over the interval $[0, 3]$. With a dataset comprising $10^4$ data points and utilizing a single input dimension, the PR-nested \sindy\ model attempts to replicate this function.

The architecture of the model is as follows. The first layer is polynomial (P), involving monomials up to the third degree (i.e. $x$, $x^2$, and $x^3$). This is followed by a linear layer, which transforms the three features into seven, with a bias term included. Subsequently, the radial layer (R) applies a series of predefined functions on a node-wise basis, belonging to the following dictionary:
\begin{equation*}
    \mathcal{F} = \left\{
    x \mapsto x,
    x \mapsto x^2,
    \arctan,
    \sin,
    \cos,
    \exp,
    x \mapsto \log(|x| + 10^{-5}),
    x \mapsto \frac{1}{1+x^2}
    \right\}.
\end{equation*}
The last linear layer maps these seven transformed features down to a single output.

The hyperparameters and training settings are as follows:
\begin{itemize}
    \item Optimizer: the Adam optimization algorithm is used, with a learning rate of $10^{-3}$.
    \item Weight initialization: the P and R layers do not use any weight initialization. The linear layer, however, applies a normal distribution initialization, where the weights are drawn from a Gaussian distribution with a specified mean ($\mu$) set to 0 and a standard deviation ($\sigma$) set to 1.5.
    \item Batch: the training is performed on $10\,000$ uniformly distributed values, with a batch size of $n_\text{batch} = 1\,000$, resulting in 10 iterations per epoch.
    \item Sparsity: Lasso regularization is applied at a rate of $10^{-3}$, and the model weights undergo pruning every~$30$ epochs with a threshold set to $0.05$.
\end{itemize}
The function uses a \texttt{patience} counter to keep track of how many consecutive epochs have passed without significant improvement in the loss function. This is calculated by comparing the relative difference in loss between the current and previous epochs to a predefined threshold~(here,~$10^{-2}$). If the relative difference is lower than this threshold for $50$ consecutive epochs ($\texttt{patience} = 50$), the training stops.
In this case, the training stopped after $300$ epochs
out of a maximum of $1\,000$.

\begin{figure}[htbp]
    \centering
    \begin{subfigure}[b]{0.49\textwidth}
        \includegraphics[width=\textwidth]{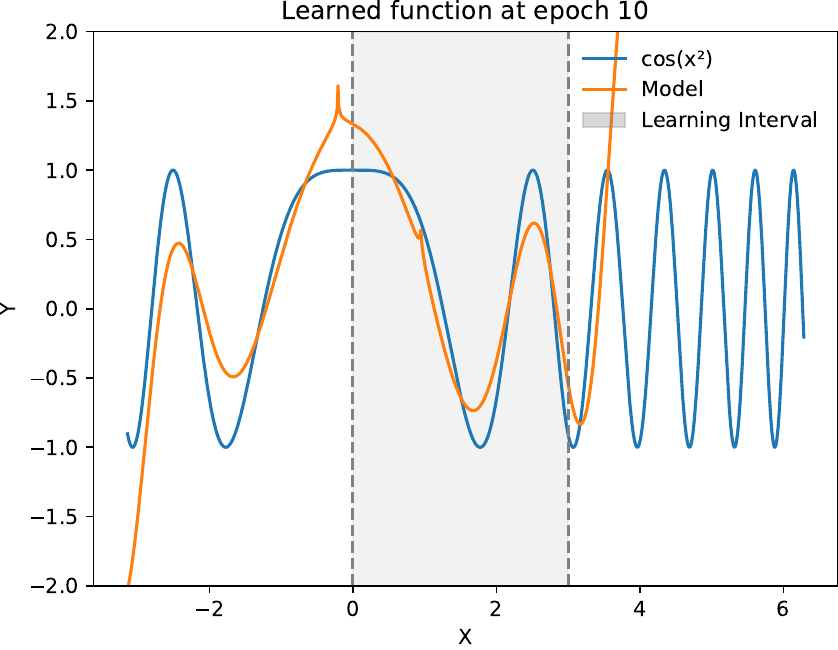}
        \caption{Epoch 10: early model predictions vs. target function, illustrating the initial learning phase.}
        \label{fig:pr_50}
    \end{subfigure}
    \begin{subfigure}[b]{0.49\textwidth}
        \includegraphics[width=\textwidth]{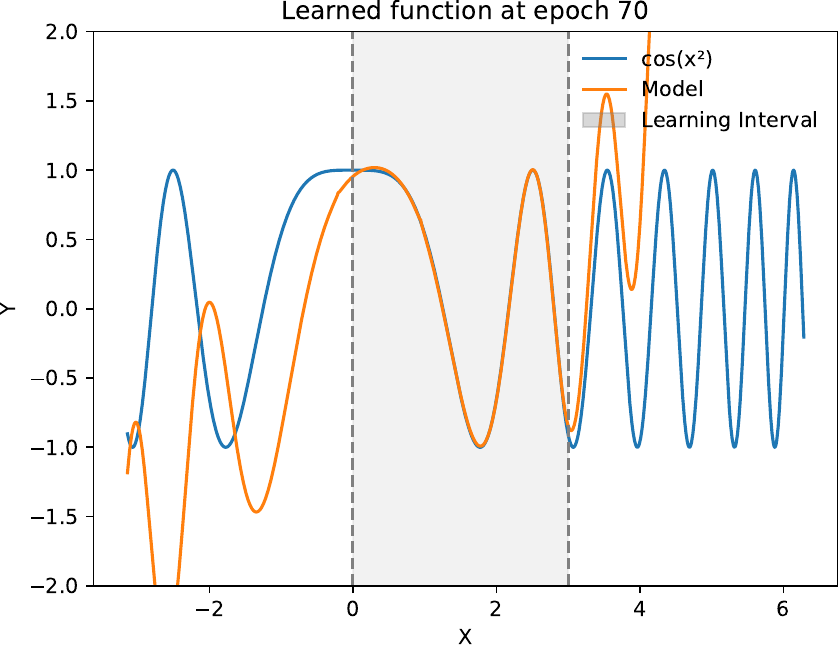}
        \caption{Epoch 70: model predictions show improved alignment with the target function as learning progresses.}
        \label{fig:pr_100}
    \end{subfigure}

    \bigskip

    \begin{subfigure}[b]{0.49\textwidth}
        \includegraphics[width=\textwidth]{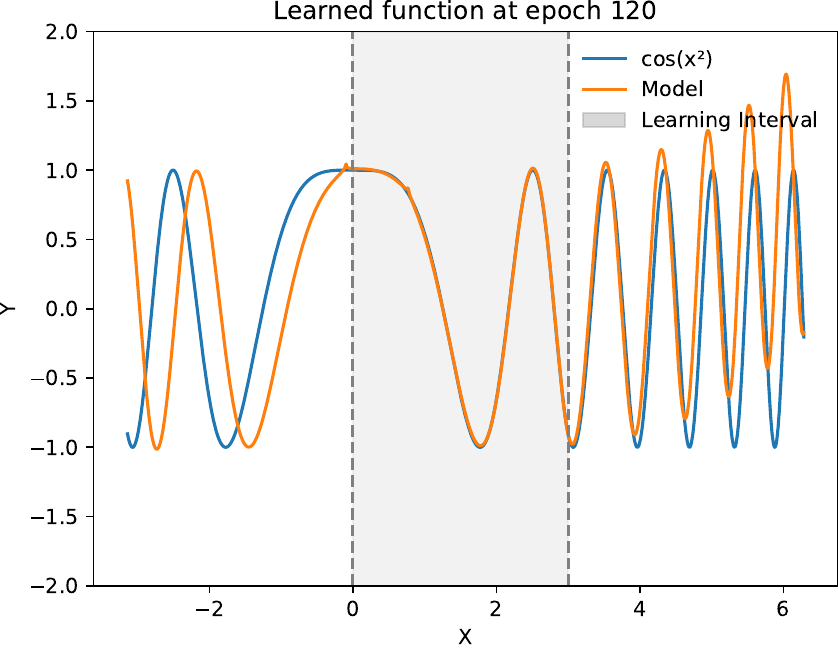}
        \caption{Epoch 120: further refined predictions, with the model beginning to capture the function's periodicity.}
        \label{fig:pr_150}
    \end{subfigure}
    \begin{subfigure}[b]{0.49\textwidth}
        \includegraphics[width=\textwidth]{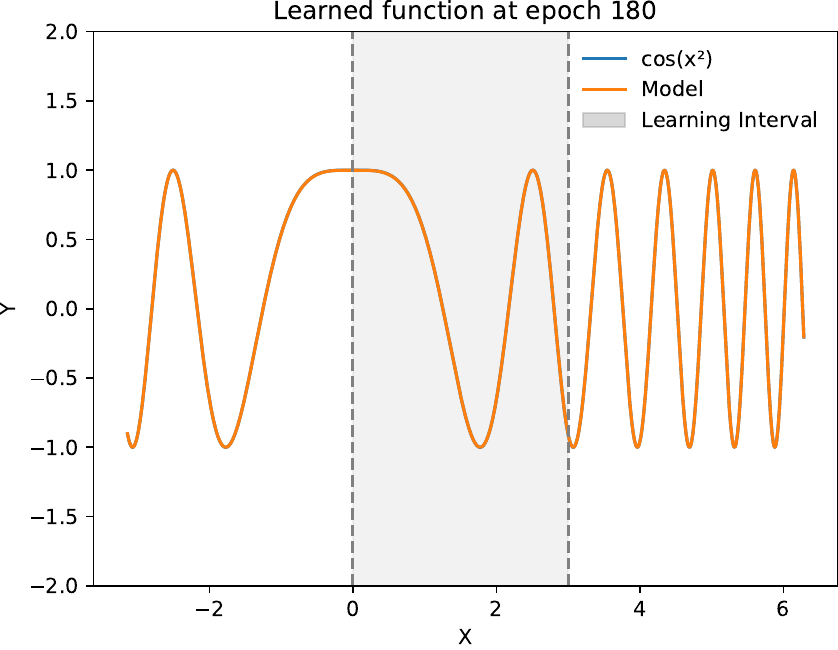}
        \caption{Epoch 180: the model has approximated the target function.}
        \label{fig:pr_250}
    \end{subfigure}
    \caption[PR Model]{Case 1 from \cref{sec:function_case_1}: evolution of the learned function $x \mapsto \cos(x^2)$ over successive training epochs. Each subfigure represents the PR model's predictions (in blue) against the target function (in orange) at epochs 10, 70, 120, and 180, showcasing the model's progressive learning and convergence towards approximation of the target function. We observe that the first epochs learn on the interval $[0, 3]$, while subsequent epochs are able to generalize.}\label{fig:pr}
\end{figure}

The learning process of the PR model with the specified parameters successfully approximated the target function with the expression $x \mapsto - \sin(x^2-1.57)$ that captures the underlying behaviour of the target function within the considered domain, see \cref{fig:pr}.
Indeed, at the end of the training,
the absolute and relative MSEs are respectively given by
\begin{equation*}
    \int_0^3 \left( -\sin (x^2 - 1.57) + \cos(x^2) \right)^2 \mathrm{d}x= 8.31 \times 10^{-7}
    \text{\quad and \quad}
    \frac{\int_0^3 \left( \sin (x^2 - 1.57) + \cos(x^2) \right)^2\mathrm{d}x}{\int_0^3 cos^2(x^2)\mathrm{d}x} = 4.92\times 10^{-7}.
\end{equation*}
We clearly see that the behaviour of the function is captured up to
(single-precision) machine error.
The sparsity of the solution is significant, with only 3 nonzero coefficients out of the original 29 parameters, underscoring the effectiveness of the model in identifying the essential structure of the function within the given domain.

\subsubsection{Case 2: trigonometric function involving
    composition using the PRP block}
\label{sec:function_case_2}

To further illustrate the capabilities of the introduced PRP-Nested \sindy\ method \eqref{eq:prp_model} described in \cref{sec:PRP_model}, we examine the same function as in \cref{sec:function_case_1}, namely $f: x \mapsto \cos(x^2)$, considered over the interval $[0, 3]$. This time, employing a dataset comprising $1\,000$ data points, the PRP-Nested \sindy\ model is tasked with replicating this function.

The model's architecture begins with a polynomial (P) layer that computes monomials up to the second degree (i.e. with monomials $x$ and $x^2$). A subsequent linear layer then transforms these two polynomial features into nine outputs, incorporating a bias term in the process. The radial (R) layer follows, applying a suite of predefined functions to each node, belonging to the following dictionary:
\begin{equation*}
    \mathcal{F} = \left\{
    x \mapsto x,
    x \mapsto x^2,
    \arctan,
    \sin, \cos,
    \exp,
    x \mapsto \sqrt{x},
    x \mapsto e^{-x^2},
    x \mapsto \log(1 + e^x)
    \right\}.
\end{equation*}
The architecture then introduces an intermediate linear layer with output size set to $k=2$, a second polynomial (P) layer still with $d=2$, i.e. with monomials $x$, $y$, $xy$, $x^2$, $y^2$, and a final linear layer, thus completing the structure of the model.

The hyperparameters match those of \cref{sec:function_case_1}, except that the Adam learning rate is set to $10^{-4}$.
We note that the loss function stagnates after around $300$ epochs.

\begin{figure}
    \centering
    \includegraphics[width=0.75\textwidth]{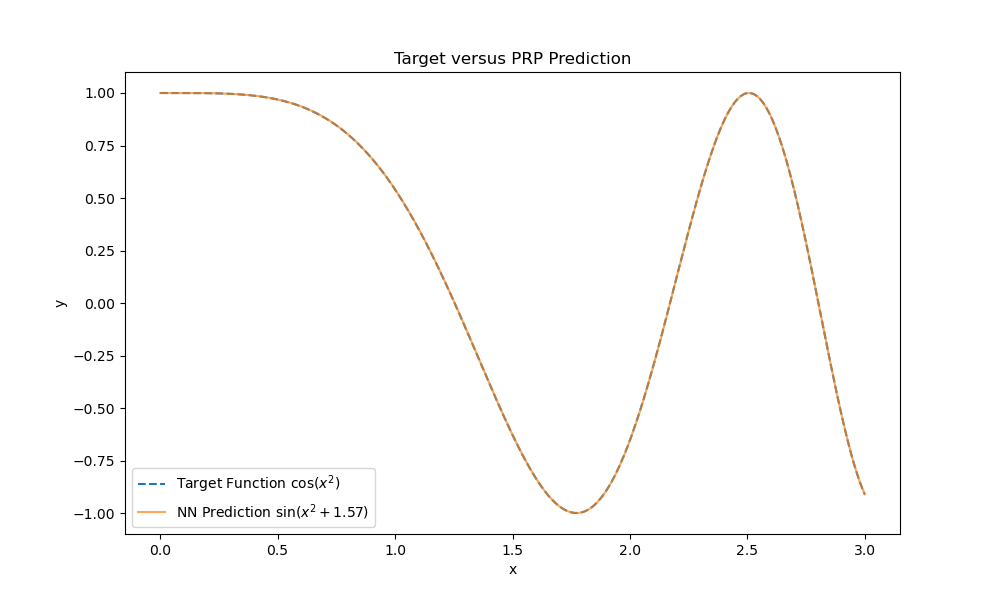}
    \caption{Case 2 from \cref{sec:function_case_2}: comparison of the results from the PRP-Nested \sindy, which yielded the function $x \mapsto \sin(x^2+1.57)$, with the target function $x \mapsto \cos(x^2)$, over the interval [0, 3].}
    \label{fig:prp}
\end{figure}

The PRP model's learning process
successfully approximated the expression $x \mapsto \sin(x^2+1.57)$. This expression effectively captures the underlying behaviour within the specified domain (\cref{fig:prp}).
After training,
the absolute and relative MSEs are respectively given by
\begin{equation*}
    \int_0^3 \left(  \cos(x^2) - \sin (x^2 + 1.57) \right)^2 \mathrm{d}x= 8.31 \times 10^{-7}
    \text{\quad and \quad}
    \frac{\int_0^3 \left(  \cos(x^2) -  \sin (x^2 + 1.57) \right)^2\mathrm{d}x}{\int_0^3 \cos^2(x^2)\mathrm{d}x} = 4.92\times 10^{-7}.
\end{equation*}
Just like before, the behaviour of the function is captured up to
(single-precision) machine error.
The resulting model exhibits notable sparsity, with only 3 nonzero coefficients from the initial 31 parameters.

\subsubsection{Case 3: Trigonometric Function Multiplication using the PRP Block}
\label{sec:example_prp_2D}

The next test case consists in learning the two-dimensional function $(x,y) \mapsto 2\sin(x)\cos(y)$ over the space domain $[-2, 2]^2$, using PRP blocks.
This example is particularly interesting as it highlights the model's ability to learn complex functions involving the multiplication of simple functions, whose product is not included in the dictionary.
This is a clear advantage over the naive \sindy\ approach.
Indeed, when two-dimensional input variables are considered,
the \sindy\ function dictionary needs to be much larger
(e.g. including all possible products of the functions in the dictionary).

Initially, a Polynomial (P) layer processes the input, generating monomials up to the second degree (i.e. $x$, $y$, $xy$, $x^2$, $y^2$). Following this, the first linear layer transforms these polynomial features from $5$ inputs to $9$ outputs, including a bias term. Next comes the Radial (R) layer, which applies a variety of nonlinear transformations to each node, belonging to the following dictionary:
\begin{equation*}
    \mathcal{F} =
    \begin{Bmatrix}
        x \mapsto \sqrt{|x|+10^{-5}},                 &
        x \mapsto x,                                  &
        \sin, \vphantom{\dfrac{1}{1+x^2}}             &
        \cos,                                         &
        \tanh,                                          \\
        \exp,                                         &
        x \mapsto \dfrac{1}{1+x^2},                   &
        x \mapsto \log(|x|+10^{-5}),                  &
        x \mapsto \exp \left(\dfrac{1}{1+x^2}\right), &
        x \mapsto \log(1 + e^x)
    \end{Bmatrix}.
\end{equation*}
After the R layer, a second P layer is applied, also containing monomials up to degree 2.
This additional polynomial processing adds depth to the feature extraction. The architecture then continues with another linear layer, which consolidates the outputs of the preceding layers. This final layer maps the 5-dimensional output from the second P layer to a single output dimension, providing the final model prediction.

The training parameters for learning the function $(x,y) \mapsto 2\sin(x)\cos(y)$ are specifically configured as follows:
\begin{itemize}
    \item Optimizer: the Adam optimization algorithm is employed with a base learning rate of $10^{-4}$, and is run for $1\,000$ epochs.
    \item Weight initialization: in the uniform weight initialization for the linear layers, the standard deviation is set to $0.5$, and the mean is maintained at $0$.
    \item Batch: the training is performed on $10\,000$ uniformly distributed values, with a batch size of $n_\text{batch} = 1\,000$, resulting in 10 iterations per epoch.
    \item Sparsity: the Lasso regularization rate is set to $0.1$, and the pruning of model weights occurs every $50$ epochs with an initial pruning at $10$ epochs, using a pruning threshold of $0.01$.
\end{itemize}

\begin{figure}[htb!]
    \centering
    \begin{subfigure}[b]{0.49\textwidth}
        \includegraphics[scale=0.35]{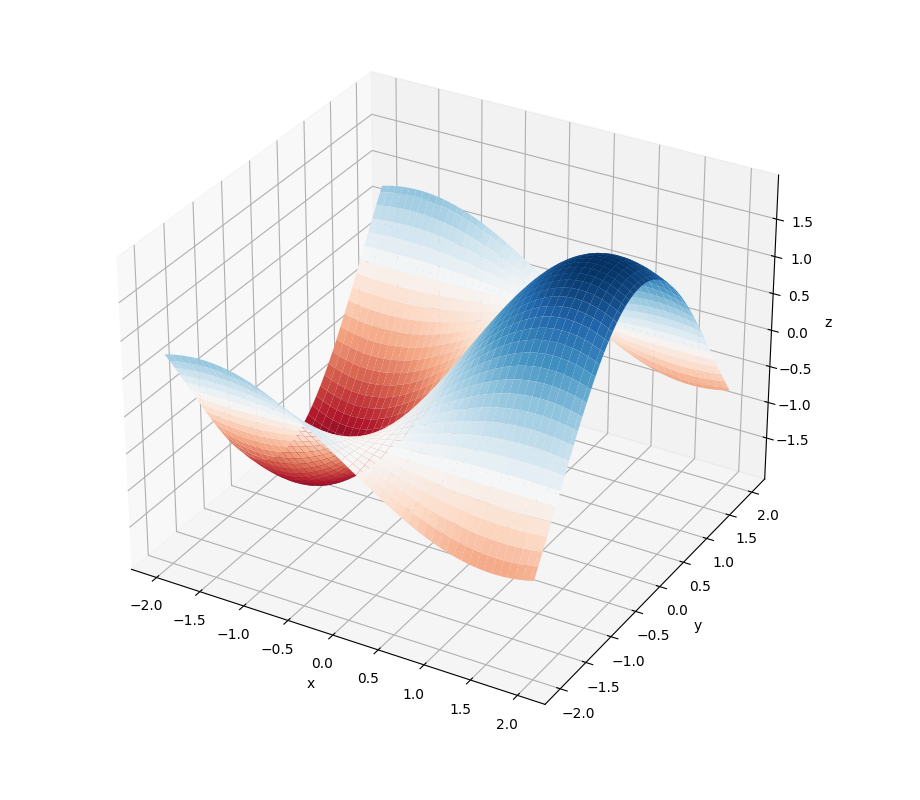}
        \caption{True model ($2\sin(x)\cos(y)$)}
        \label{fig:prp_true}
    \end{subfigure}
    \begin{subfigure}[b]{0.49\textwidth}
        \includegraphics[scale=0.35]{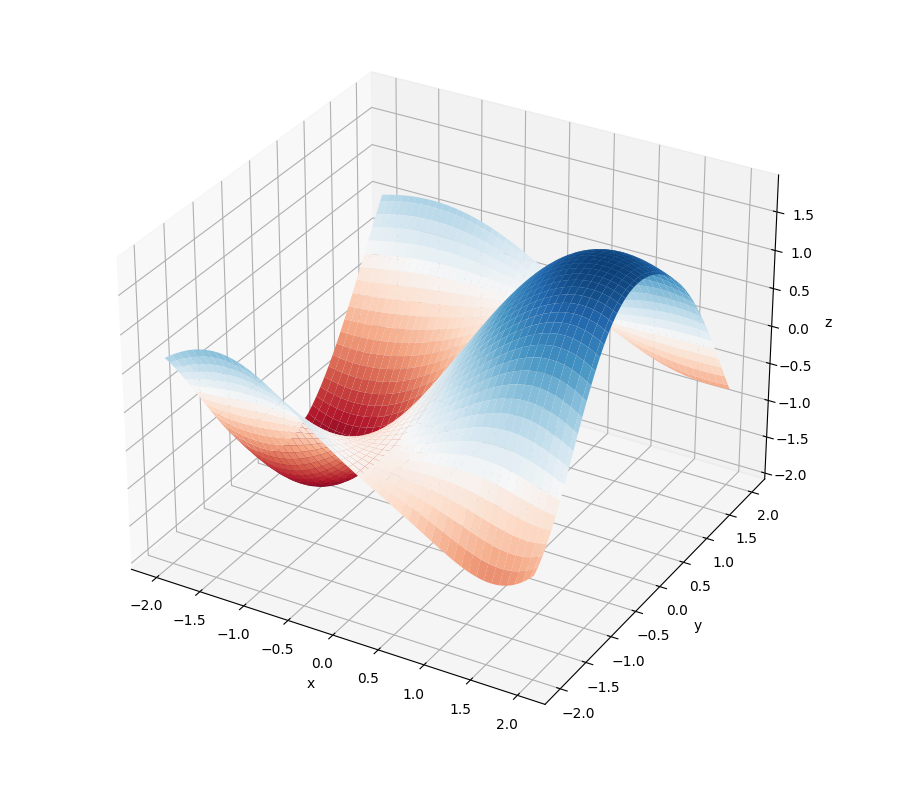}
        \caption{Predicted model}
        \label{fig:prp_pred}
    \end{subfigure}
    \caption[PRP Model]{Case 3: comparison of the true and predicted models. A visual representation of the target function $(x,y) \mapsto 2\sin(x)\cos(y)$ and its learned approximation by the PRP model, given by \eqref{eq:prp_2D_result}, over the specified domain.}
    \label{fig:prp_2D}
\end{figure}

Utilizing these parameters, the PRP blocks successfully learn the following function:
\begin{equation}
    \label{eq:prp_2D_result}
    \begin{aligned}
        (x,y) \mapsto & -0.77 \left(1 - 0.612 \sin \left(0.032x^2 - 0.162xy + 0.998x + 0.035y^2 - y - 0.14\right)\right)^2 \\
                      & + 2.01 \left(\cos \left(0.039xy - 0.499x - 0.497y + 0.809\right)\right)^2
    \end{aligned}
\end{equation}
with sparsity accounted for by 14 nonzero coefficients out of the original 32 parameters.
A comparison of this function with the target function is illustrated in \cref{fig:prp_2D}.
The MSE between the target function and its approximation is
$5.69 \times 10^{-2}$.

\subsubsection{Case 4: Perimeter of an ellipse}

Calculating the perimeter of an ellipse is a well-studied topic for which the solution cannot be expressed in terms of elementary functions.
Let us first define an ellipse as the set of points $(x, y)$ such that:
\begin{equation*}
    \label{eq:ellipse}
    \frac{x^2}{a^2} + \frac{y^2}{b^2} = 1,
\end{equation*}
which can be described by the parametric equations:
\begin{equation*}
    \label{eq:parametric}
    \begin{cases}
        x = a \cos(\alpha), \\
        y = b \sin(\alpha),
    \end{cases}
\end{equation*}
with $\alpha\in[0, 2\pi)$.
The perimeter of an ellipse can then be expressed as:
\begin{equation}
    \label{eq:perimeter}
    P(a, b) = 4 \int_0^{\frac{\pi}{2}} \sqrt{a^2 \cos^2(\alpha) + b^2 \sin^2(\alpha)} \, d\alpha,
\end{equation}
because the differential arc length of our parametric equation is $ds = -\sqrt{((-a \sin(\alpha))^2 + (b \cos(\alpha))^2)} d\alpha$ and the four quadrants have the same length.
Several approximations are known to approach the perimeter of an ellipse, such as the one proposed by Ramanujan:
\begin{equation}
    \label{eq:ramanujan}
    P(a, b) \approx \pi (3(a + b) - \sqrt{(3a + b)(a + 3b)}).
\end{equation}

Since any rescaled ellipse remains an ellipse, we assume without loss of generality that $b = 1$ for the rest of this experiment.
To assess the relevance of our Nested \sindy\ technique, we will compare its performance against the two reference solutions: Ramanujan's approximation, and a linear interpolation between the two endpoints of the interval $[1, 30]$.
Indeed, using a linear interpolation makes sense for large values of $a$, since $P$ is equivalent to $4a$ when $a \to \infty$, as can be seen by inspecting \eqref{eq:perimeter}.

To find formulas that approximate the perimeter of an ellipse, we train two PRP models.
The first one is trained on the interval $[1, 5]$ and the second one on the interval $[1, 25]$.
We run 50 learning sessions for each model and keep the model with the shortest associated formula.
The only difference between the learning sessions is the seed used for initializing the random number generator.
A learning session consists of two trainings, each with its own parameters:
\begin{itemize}
    \item Optimizer: the Adam optimization algorithm is used for both trainings. The first training is run for $100$ epochs with a base learning rate of $10^{-1}$, while the second one includes $1\,500$ epochs and a learning rate of $10^{-3}$.
    \item Weight initialization: the initial weights are uniformly sampled with mean $0$ and standard deviation $0.5$.
    \item Batch: the training is performed on $1\,000$ uniformly distributed values, using a single batch.
    \item Sparsity: the Lasso regularization rate is set to $10^{-1}$ in the first training and $10^{-2}$ in the second one. The pruning occurs after $30$ epochs below the threshold of $5\times 10^{-2}$ in both trainings.
\end{itemize}
Overall, the first training aims to move faster near a solution, while the second one aims to converge to a precise solution.

In its best session, the first model found a basic quadratic polynomial:
\begin{equation}
    \label{eq:quadratic}
    P_\text{quadratic}(a) = 0.061(a + 0.544)^2 + 3.28a + 2.72.
\end{equation}
However, other relatively sparse solutions were found,
such as a model with 9 nonzero parameters
\begin{equation*}
    \label{eq:other_1}
    P_1(a) = 1.65a + 0.553(0.485a + 0.135\log(0.817|a^2|) + 1)^2 + 0.459\log(0.817|a^2|) + 3.4,
\end{equation*}
or one with 15 nonzero parameters
\begin{equation*}
    \label{eq:other_2}
    \begin{aligned}
        P_2(a) & = 0.535a + 0.966(0.394a + 0.721\arctan(0.278a^2 + 0.393) + 1                             \\
               & + 0.111\exp(-0.063a^4))^2 + 0.978\arctan(0.278a^2 + 0.393) + 1.36 + 0.15\exp(-0.063a^4),
    \end{aligned}
\end{equation*}
thus demonstrating the ability of the model to find a variety of approximations.
The MSE obtained at the end of the training are $2.26 \times 10^{-3}$, $2.69 \times 10^{-3}$, and $3.30 \times 10^{-3}$ for the quadratic, $P_1$, and $P_2$ models, respectively.

The second model found a more complex expression which is too long to be displayed here.
In this second model, 25 out of the 64 parameters are nonzero, corresponding to approximately 40\% sparsity.
The MSE obtained at the end of the training is $1.97$.
This metric, as the previous ones, is biased because the training is stopped at an arbitrary step and there is no final tuning step.
Hence, we perform a final tuning step, by writing the final expression of the model and performing a gradient descent to only minimize the MSE.
The later mentioned results include this final tuning step and show that the effective MSE can be significantly reduced.

\begin{figure}[htbp]
    \centering
    \includegraphics[width=0.75\textwidth]{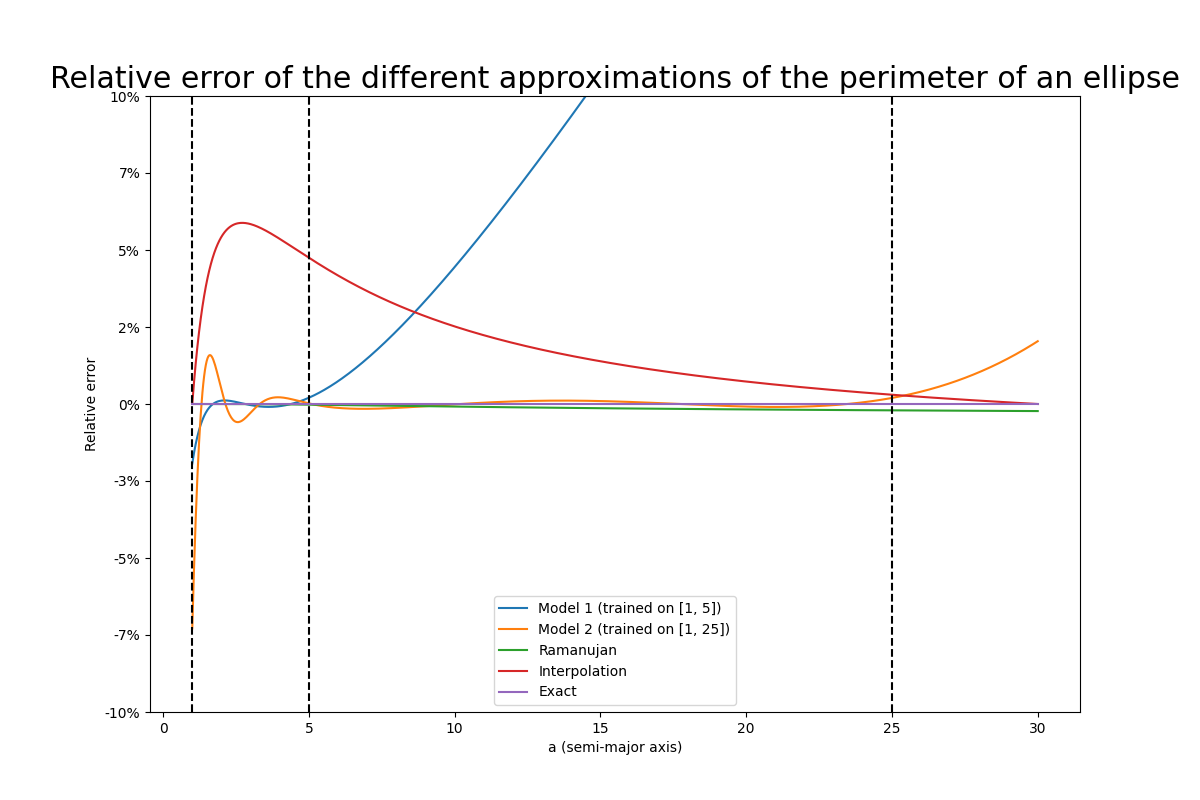}
    \caption{Case 4: relative error of the different approximations of the ellipse perimenter on the intervals $[1, 5]$, $[1, 25]$ and $[1, 30]$.}
    \label{fig:ellipse}
\end{figure}

\begin{table}[htbp]
    \caption {Case 4: mean squared error of the different approximation of the ellipse perimeter on the intervals $[1, 5]$, $[1, 25]$ and $[1, 30]$. The closest approximations are denoted in bold.} \label{tab:ellipse}
    \centering
    \begin{tabular}{rccc}
        \toprule
                      & $[1, 5]$
                      & $[1, 25]$
                      & $[1, 30]$                      \\
        \cmidrule(lr){1-4}
        Model 1       & $7.24 \times 10^{-4}$
                      & $1.05 \times 10^{2}$
                      & $2.61 \times 10^{2}$           \\
        Model 2       & $7.43 \times 10^{-3}$
                      & $\mathbf{3.66 \times 10^{-3}}$
                      & $2.73 \times 10^{-1}$          \\
        Ramanujan     & $\mathbf{2.78 \times 10^{-6}}$
                      & $9.84 \times 10^{-3}$
                      & $\mathbf{1.82 \times 10^{-2}}$ \\
        Interpolation & $4.70 \times 10^{-2}$
                      & $4.75 \times 10^{-1}$
                      & $5.46 \times 10^{-1}$          \\
        \bottomrule
    \end{tabular}
\end{table}

We display the relative error of each approximation on different intervals in \cref{fig:ellipse}.
\cref{tab:ellipse} reports the MSE of each approximation on the intervals $[1, 5]$, $[1, 25]$ and $[1, 30]$.
\textit{Model 1} and \textit{model 2} correspond to our two models (with \textit{model 1} corresponding to \eqref{eq:quadratic}), \textit{Ramanujan} corresponds to the approximation proposed by Ramanujan and \textit{Interpolation} corresponds to the linear interpolation between the two endpoints of the interval $[1, 30]$.
The two models perform well on the interval they were trained on.
The first model performs better than the linear interpolation on $[1, 5]$, but worse than the Ramanujan approximation.
The second model outperforms both the linear interpolation and the Ramanujan approximation on $[1, 25]$, but starts diverging on $[1, 30]$.

Overall, this experiment provides an insight into the ability of the Nested \sindy\ approach to discover an approximation on a classical problem without knowledge other than the data points.
It demonstrates that this method can be used to discover an approximation of a complex function with little effort.
The main downside of this approach is the lack of guarantee to find a sparse solution.
The initial choice of the model's coefficients appears to have a significant impact on the final solution.

\subsection{ODE discovery}
\label{subsec:ode_discovery}

\begin{figure}[htbp]
    \vspace{-3cm}
    \begin{subfigure}[b]{0.49\textwidth}
        \centering
        \tiny{
            \centering
            \def\svgwidth{9cm}
            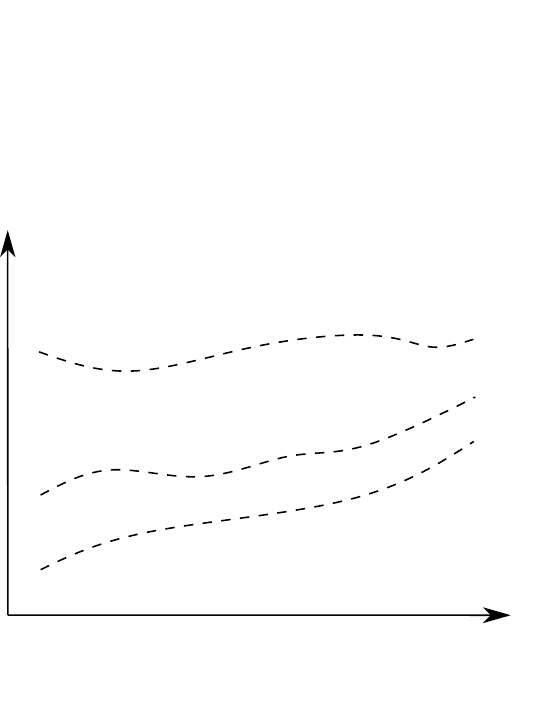}
        \vspace{-1.65cm}
        \caption{Illustration of the full data set for ODE discovery. \\ ~ }
        \label{fig:datasetode}
    \end{subfigure}
    \begin{subfigure}[b]{0.49\textwidth}
        \centering
        \tiny{
            \centering
            \def\svgwidth{9cm}
            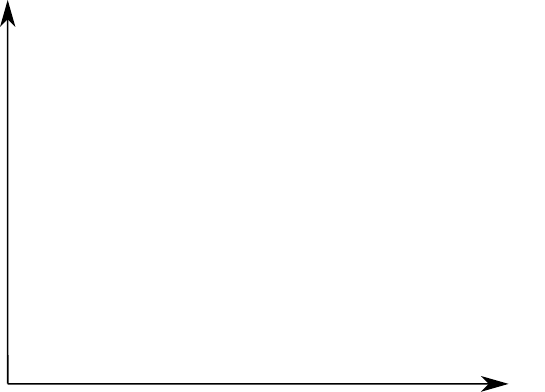
        }
        \centering
        \caption{Reduced dataset for one training step, with $n_\text{batch} = 2$ and $l_\text{batch} = 4$ (batches are highlighted in blue).}
        \label{fig:batchesode}
    \end{subfigure}
    \caption{Illustration of data set and batches for the training of the Neural Nested SINDy model for ODE discovery.}
    \label{fig:data_ode}
\end{figure}

In this section, we extend the proposed approach to try to discover an unknown autonomous ODE
\begin{equation}
    \label{eq:unknown_ODE}
    x'(t) = f(x(t)),
\end{equation}
based on time data associated to $K$ discrete trajectories in time.
The full dataset (see \cref{fig:datasetode}) is given by
\begin{equation*}
    \mathcal{X}=\bigcup_{i=1}^K \mathcal{X}^k,
\end{equation*}
where we have defined the $k$\textsuperscript{th} trajectory as follows:
\begin{equation*}
    \mathcal{X}^k=\left\{x_{1}^{(k)},\dots,x_{N}^{(k)}\right\},
\end{equation*}
and the subscripts indicate the time steps.
The goal is to provide an approximation $f_{\theta}$ of~$f$ in \eqref{eq:unknown_ODE}.

To approximate the dynamical system, we will adapt the method developed by the authors of \cite{lee2021structurepreserving}.
They suggest a coupling between \sindy\ and Neural ODE approaches.
Among other things, this method appears to be robust even when the dataset is collected at large or irregular time steps, or contains additive noise.

The idea is, first, to minimize
\begin{equation*}
    L(\theta)= \frac{1}{K}\sum\limits_{k=1}^{K}\sum\limits_{i=1}^{N} \left\lVert x_{\theta,i}^{(k)} - x_i^{(k)} \right\rVert_2 + \lambda \left\lVert  \theta\right\rVert_1
    \label{eq: lossfunctionODE}
\end{equation*}
where $(x_{\theta,i}^{(k)})_{i \ge 0}$ is the solution of the following Cauchy problem:
\begin{equation}
    \begin{cases}
        x'(t) = f_{\theta}(x(t)), \\
        x(0) = x_0^{(k)},
    \end{cases}
    \label{eq: edo avec ftheta}
\end{equation}
where $f_{\theta}$ is given by the Nested \sindy\ network. In practice, the ODE \eqref{eq: edo avec ftheta} is solved numerically using an ODE solver (a fifth-order Runge-Kutta method in our application). Hence, $(x_{\theta,i}^{(k)})_{i \ge 0}$ is an \textit{approximate} solution of \eqref{eq: edo avec ftheta}.
To train the model, the strategies described in \cref{sec:training} are used.
One must be careful about how the dataset is used to train the model. To train efficiently the Nested \sindy\ network, we divide the dataset in batches at each training step. To do that, we randomly sample $n_\text{batch}$ trajectories from the dataset $\mathcal{X}$. Then, we randomly sample initial points from the selected $n_\text{batch}$ trajectories, in such a way that the reduced dataset consists of $n_\text{batch}$ trajectories of fixed lengths $l_\text{batch}$, as shown in \cref{fig:batchesode}. The model parameters are then updated by minimizing the following reduced loss function:
\begin{equation*}
    L_\text{red}(\theta) =
    \frac{1}{n_\text{batch}}
    \sum_{k=1}^{n_\text{batch}}
    \sum_{i=i_0^{(\sigma(k))}}^{i_0^{(\sigma(k))}+l_\text{batch}-1}
    \left\lVert x_{\theta,i}^{(\sigma(k))} - x_i^{(\sigma(k))} \right\lVert_2 +
    \lambda \left\lVert  \theta\right\lVert_1
    \label{eq: reducedlossfunctionODE}
\end{equation*}
where $
    \mathcal{X}^{\sigma(1)},
    \dots,
    \mathcal{X}^{\sigma(k)},
    \dots,
    \mathcal{X}^{\sigma (n_\text{batch})}
$ are the sampled trajectories, with the associated sampled initial points $
    i_0^{(\sigma(1))},
    \dots,
    i_0^{(\sigma(k))},
    \dots,
    i_0^{(\sigma(n_\text{batch}))}
$.

Through the applications, we will attempt to demonstrate the added value of the Nested \sindy\ algorithm compared to the classical one, considering cases where the function $f$ is a composition or product of elementary functions.
It should be noted that we will only consider scalar ODEs, but it is entirely possible to explore systems of ODEs with the same method.

\subsubsection{Case 1: A trigonometric function involving composition using the PR block}
\label{sec:ode_case_1}

Firstly, we try to recover the following ODE:
\begin{equation}
    x'(t)=\sin(x^2)
    \label{eq: sinx2ODE}
\end{equation}
The function $f: x \mapsto \sin(x^2)$ is a composite function, that could be learned by a PR block.
The dataset consists of $K=100$ trajectories with $N=500$ constant time steps for each trajectory on the time interval $[0,1]$. The initial condition of each trajectory is uniformly sampled between $-3$ and $3$.

\begin{figure}[htbp]
    \centering
    \begin{subfigure}{\textwidth}
        \includegraphics[width=0.49\textwidth]{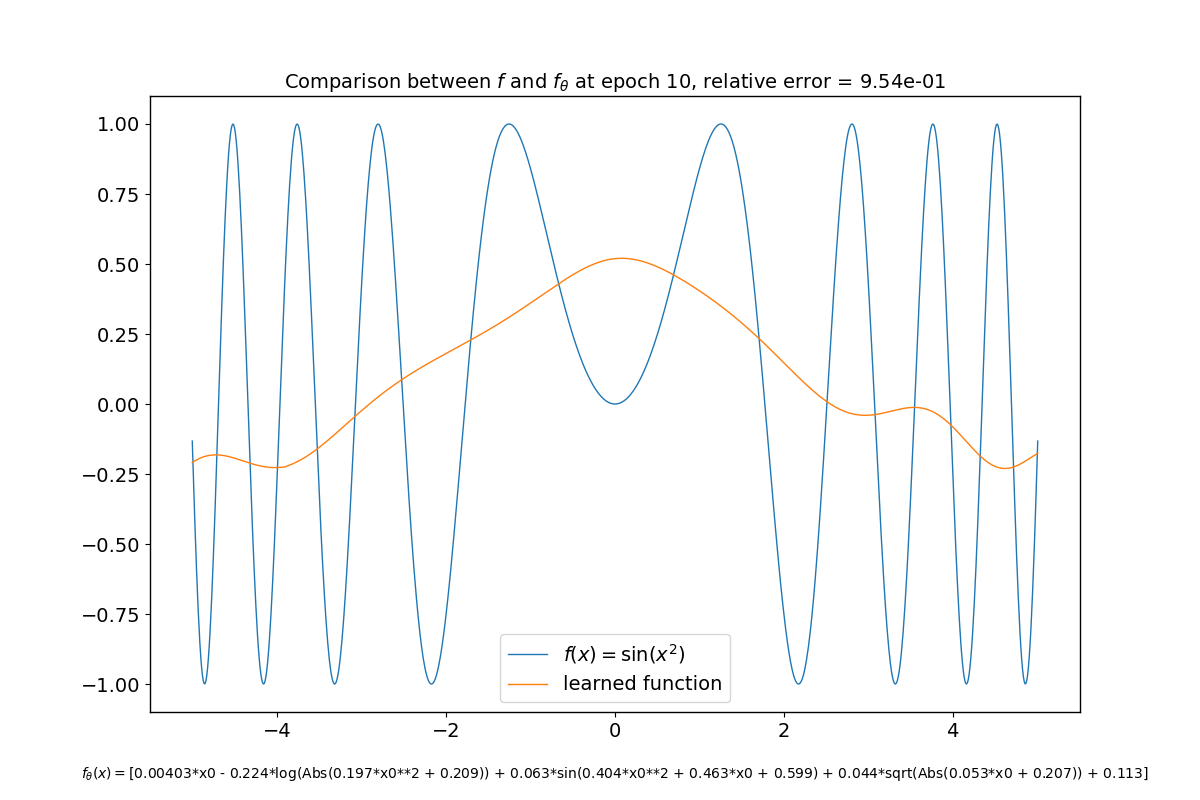}
        \hfill
        \includegraphics[width=0.49\textwidth]{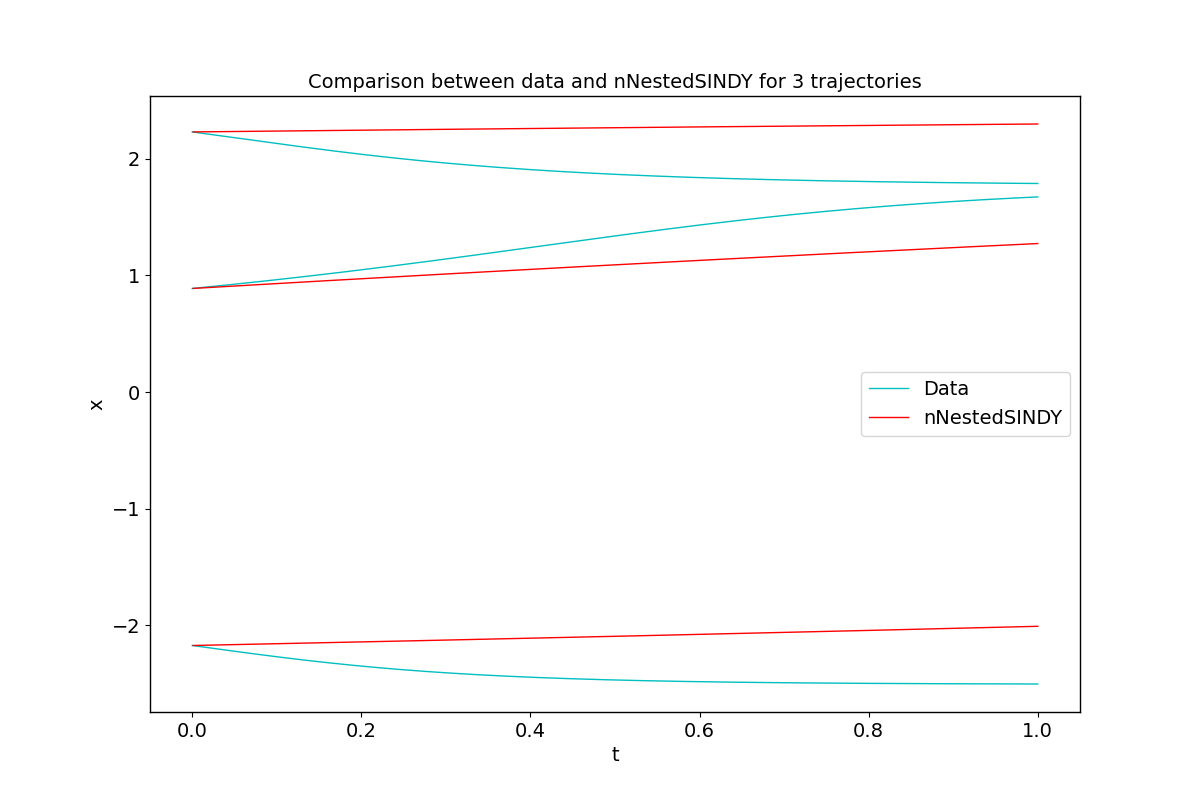}
        \caption{After 10 epochs:
            $f_\theta(x) = 0.004x -0.224\log (|0.197x^2+0.209|) + 0.063\sin(0.404x^2+0.463x+ 0.599) + 0.044 \sqrt{| 0.053x+0.207|} + 0.113$.}
    \end{subfigure}
    \medskip
    \begin{subfigure}{\textwidth}
        \includegraphics[width=0.49\textwidth]{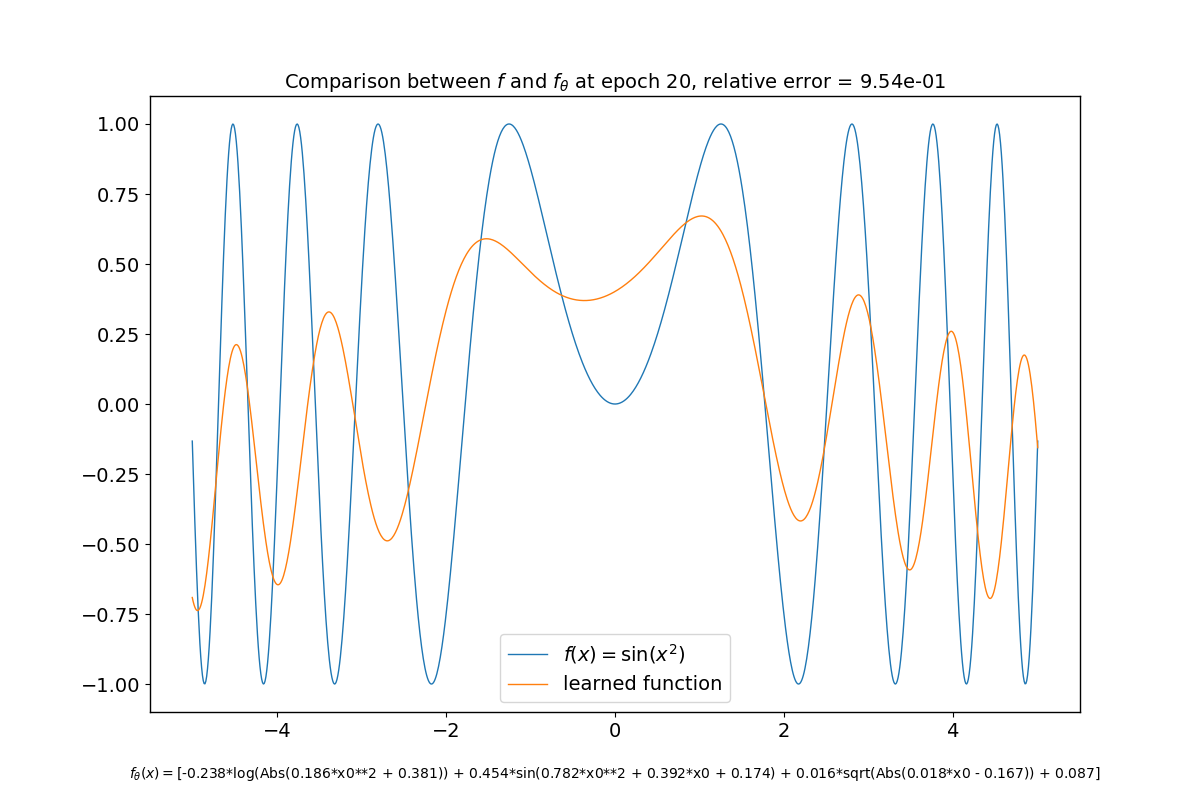}
        \hfill
        \includegraphics[width=0.49\textwidth]{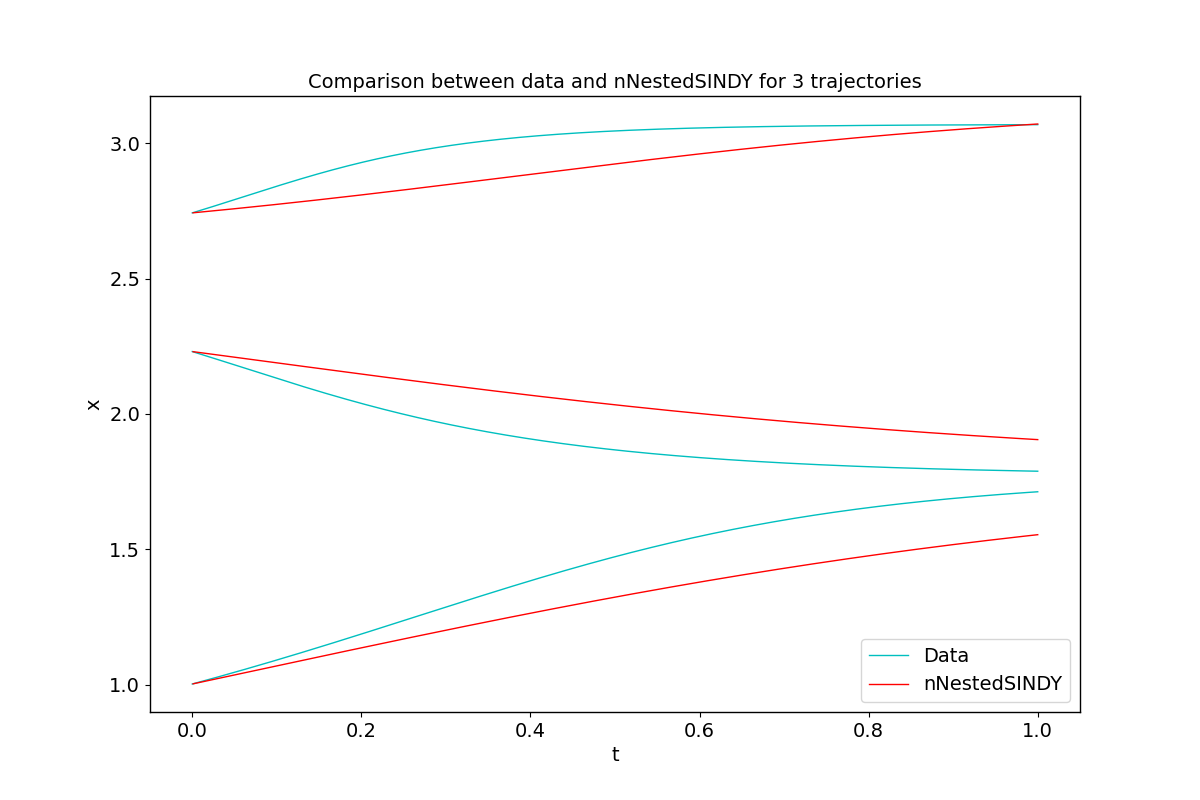}
        \caption{After 20 epochs:
            $f_\theta (x) = -0.238\log(|0.186x^2+0.381|)+0.454\sin(0.782x^2+0.392x+0.174)+0.016\sqrt{|0.018x-0.167|}+0.087$.}
    \end{subfigure}

    \caption{Case 1 from \cref{sec:ode_case_1}: training of Nested \sindy\ to recover equation  \eqref{eq: sinx2ODE}. In the left figures, we compare $f$ (blue curve) and $f_\theta$ (orange curve), while in the right figures, we compare the trajectories of dynamic systems  \eqref{eq: edo avec ftheta} (red curves) and \eqref{eq: sinx2ODE} (blue curves) for three different initial conditions. After 20 epochs, the model is still far from the true formula.}
    \label{fig: training ode sinx2_part1}
\end{figure}

\begin{figure}[htbp]
    \centering
    \begin{subfigure}{\textwidth}
        \includegraphics[width=0.49\textwidth]{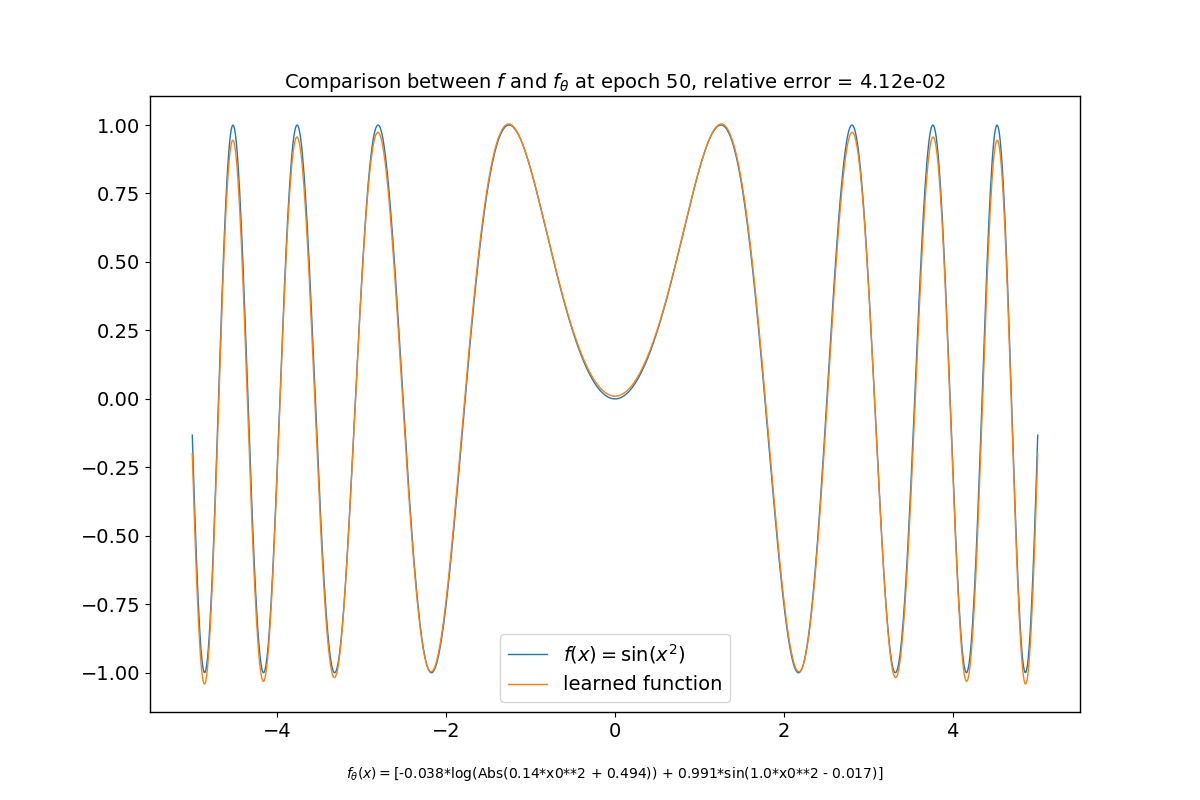}
        \hfill
        \includegraphics[width=0.49\textwidth]{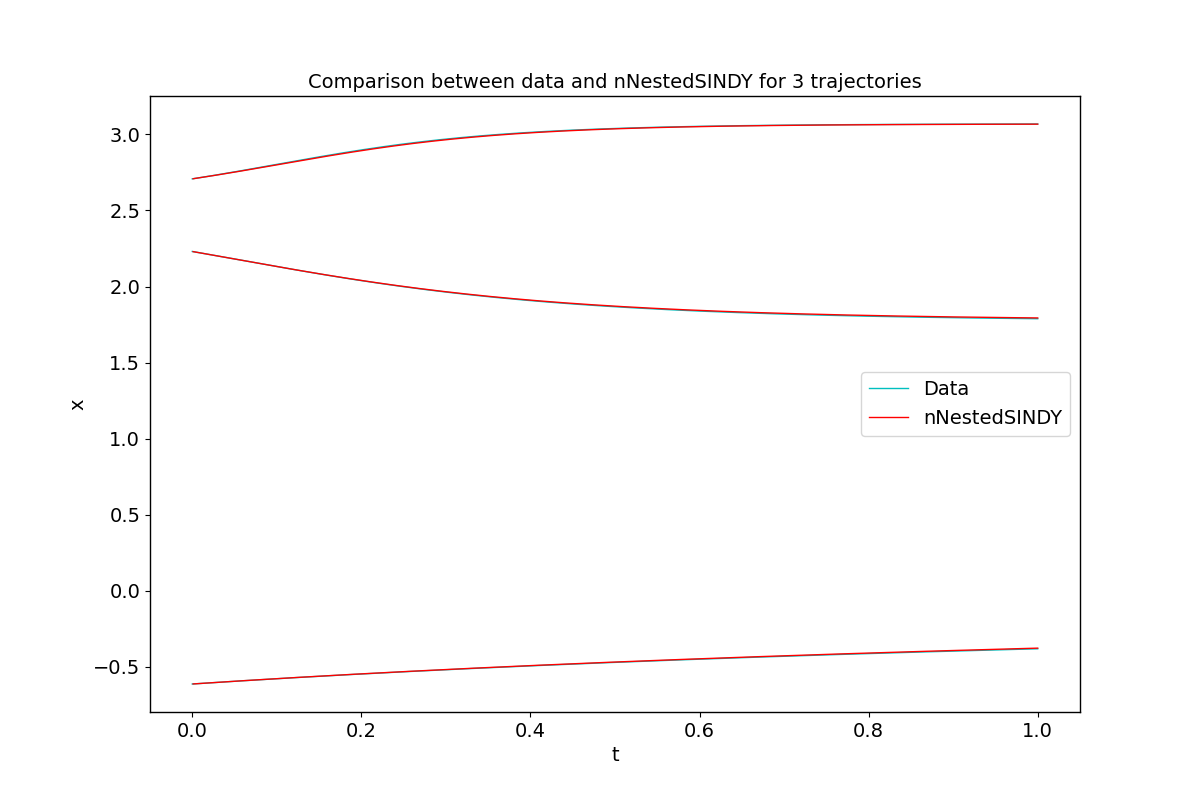}
        \caption{After 50 epochs:
            $f_\theta (x) = -0.038\log(|0.14x^2+0.494|)+0.991\sin(x^2-0.017)$.}
    \end{subfigure}
    \medskip
    \begin{subfigure}{\textwidth}
        \includegraphics[width=0.49\textwidth]{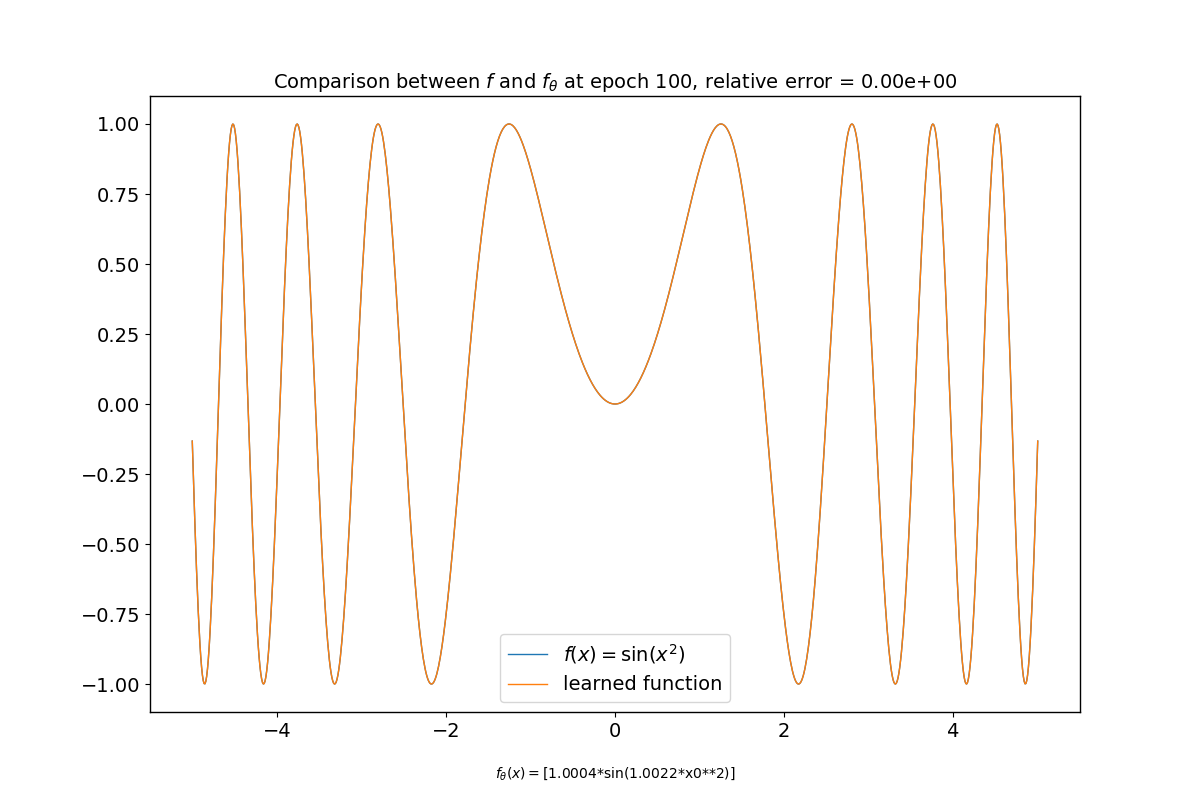}
        \hfill
        \includegraphics[width=0.49\textwidth]{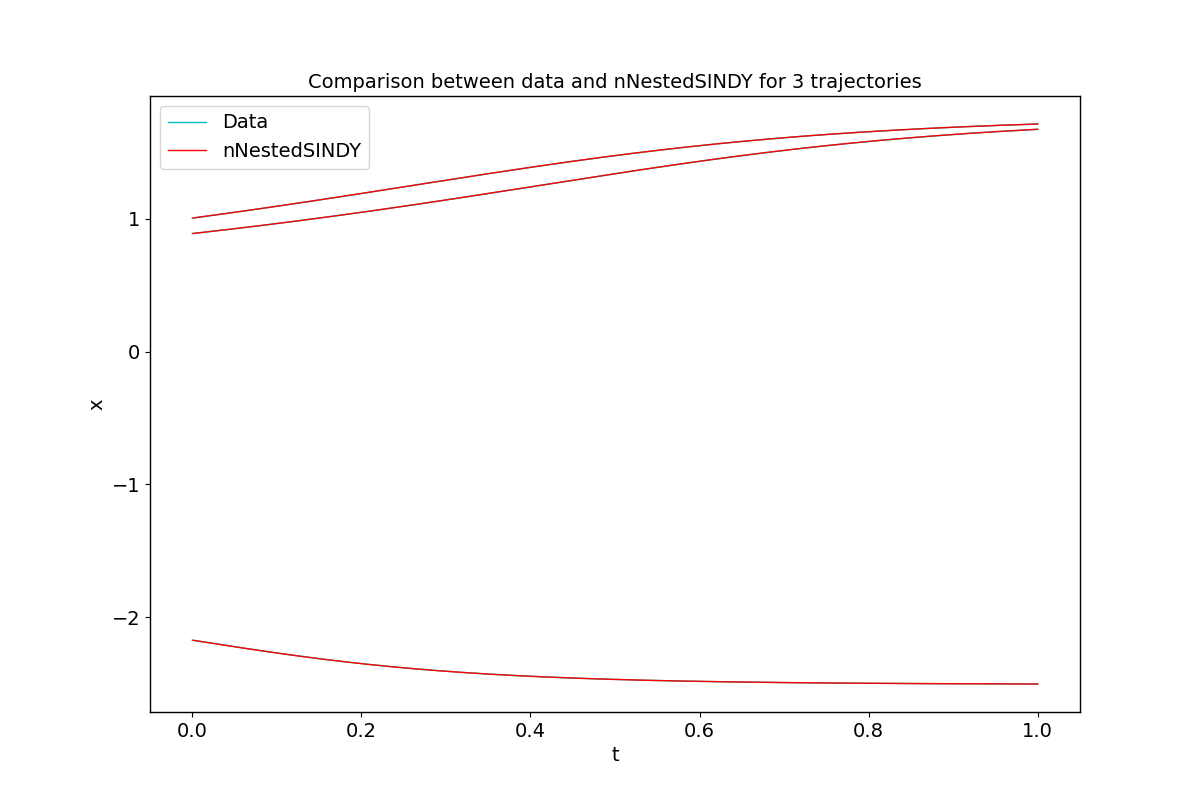}
        \caption{After 100 epochs:
            $f_\theta(x) = 1.0004\sin(1.0022x^2)$.}
    \end{subfigure}

    \caption{Case 1 from \cref{sec:ode_case_1}: training of Neural Nested SINDy to recover equation~\eqref{eq: sinx2ODE}. In the left figures, we compare $f$ (blue curve) and $f_\theta$ (orange curve), while in the right figures, we compare the trajectories of dynamic systems  \eqref{eq: edo avec ftheta} (red curves) and \eqref{eq: sinx2ODE} (blue curves) for three different initial conditions. We observe that the model approximates well the right formula after 100 epochs.}
    \label{fig: training ode sinx2_part2}
\end{figure}

Our network consists of a single PR block, where the polynomial layer contains monomials up to degree $2$, and where the radial layer contains the following functions:
\begin{equation*}
    \mathcal{F} = \{ \
    x \mapsto x, \
    x \mapsto x^2, \
    x \mapsto x^3, \
    x \mapsto \sin (x), \
    x \mapsto \log( | x |), \
    x \mapsto \sqrt{ |x| } \
    \}.
\end{equation*}
The hyperparameters and training settings are as follows:
\begin{itemize}
    \item Optimizer: we use the Adam optimizer, with learning rate set to $0.01$ and multiplied by $0.999$ every $10$ epochs. We add noise to the gradient when the loss do not decrease since $100$ epochs, as explained in \eqref{eq:add_noise_to_grad} with $\alpha = 16$.
    \item Weight initialization: all the weights are initialized to $0.2$.
    \item Batch: we set $n_\text{batch} = 30$ and $l_\text{batch} = 5$.
    \item Sparsity: the Lasso parameter is oscillating around $ \lambda_0 = 10^{-4}$ (see \eqref{eq : lasso oscillations} for details), and pruning parameters are set to $\text{MSE}_\text{prune} = 0.01$, $\varepsilon_\text{prune} = 0.01$, and $n_\text{prune}=2$.
\end{itemize}
We train the model for $100$ epochs. The results for 10 and 20 epochs are displayed in \cref{fig: training ode sinx2_part1},
where we observe that the model has not yet found a satisfactory approximation.
The model retrieves the right formula after 100 epochs, as depicted on \cref{fig: training ode sinx2_part2} and in \cref{tab:results ode sinx2}.

This example showcases the PR model's ability to accurately learn the velocity of an ODE from trajectories. To further validate this ability, future research directions include conducting tests by increasing the dictionary size and using reduced data (fewer trajectories), sparser data (lower sampling frequency in the data), and noisy data.

\begin{table}[htbp]

    \caption {Case 1 from \cref{sec:ode_case_1}: evolution of the training of Neural Nested SINDy to recover equation \eqref{eq: sinx2ODE}. The table contains the formulas of $f_ \theta$ for epochs greater than $50$, the mean squared error between $f_\theta$ and $f$, and the number of pruned parameters in the PRP block, at epochs 10, 20, 50, 100.}
    \label{tab:results ode sinx2}
    \centering
    \begin{tabular}{ccccc}
        \toprule
        Epoch
         & $f_\theta$
         & MSE on $[-3,3]$                                                                   & MSE on $[-5,5]$       & \makecell{pruned \\parameters} \\
        \cmidrule(lr){1-5}
        10
         & expression too long to display
         & $3.67 \times 10^{-1}$                                                             & $4.16 \times 10^{-1}$ & 12               \\
        20
         & expression too long to display
         & $2.91 \times 10^{-1}$                                                             & $4.16 \times 10^{-1}$ & 14               \\
        50
         & $\begin{aligned} -0.04 \log (\lvert 0.14x^2 + 0.49 \lvert) + 0.99\sin(x^2 - 0.02)
                \end{aligned}$
         & $2.29 \times 10^{-4}$                                                             & $7.72 \times 10^{-4}$ & 19               \\
        100
         & $\begin{aligned}
                    1.00 \sin(1.00x^2)
                \end{aligned}$

         & $3.54 \times 10^{-5}$                                                             & $3.00 \times 10^{-4}$ & 23               \\

        \bottomrule
    \end{tabular}
\end{table}

\subsubsection{Case 2: The Gompertz model}
\label{sec:ode_case_2}

The Gompertz model was initially designed to describe human mortality, but is now applied in various fields, especially in biology, see for instance the review paper~\cite{TjoTjo2017}. We introduce this model in the context of tumor growth. The evolution of the tumor volume $x(t)$ is governed by the following ODE:
\begin{equation}
    \label{eq: GompertzODE}
    x'(t) = r x(t)\log \left(\frac{k}{x(t)} \right) = rx(t)\log (k) - r x(t) \log (x(t)),
\end{equation}
where $k$ is the maximum size that can be reached by the tumor, and $r$ is a constant linked to the cells' proliferative capacity. In the following, we will set $k=1$ and $r=2$. By setting $k=1$, we slightly simplify the formula of the growth speed, which becomes:
\begin{equation}
    f(x) = -2x\log( x ) \,.
    \label{eq: Gompertz_speed}
\end{equation}
The exact formula of the growth speed function $f$ can only be recovered with a PRP block, to reconstruct the multiplication between the two basis functions $x \mapsto x$ and $\log$.

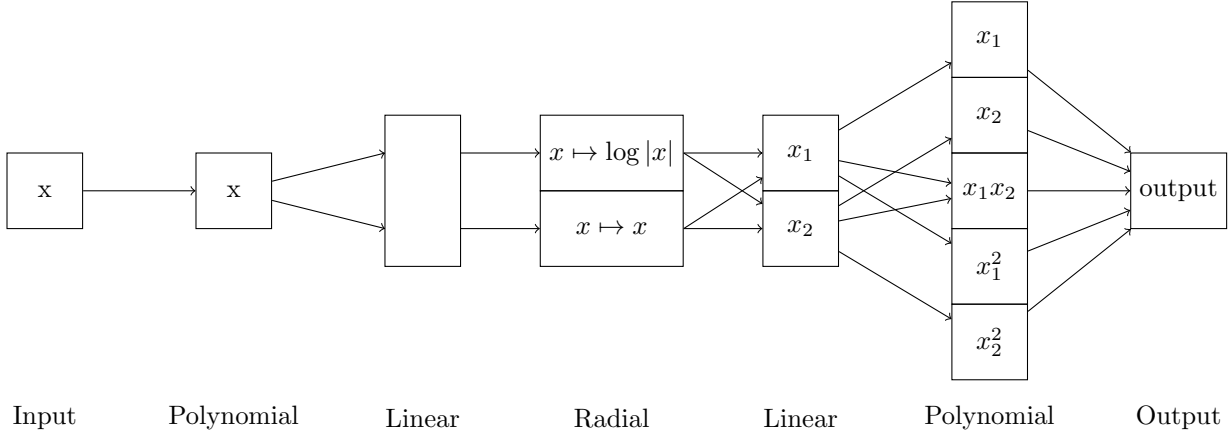
\begin{figure}[htbp]
    \centering
    \begin{tikzpicture}
        \tikzset{
            box/.style={draw, thick, rectangle},
            cell/.style={draw, rectangle, minimum height=1cm, minimum width=1cm}
        }

        \node[cell] (x_in) at (0,0cm) {x};
        \node[inner sep=0pt, outer sep=0pt, right=0cm of x_in] (input) {};

        \node[cell] (x_poly) at (2.5cm,0cm) {x};
        \node[inner sep=0pt, outer sep=0pt, right=0cm of y2_poly] (polynomial) {};

        \node[cell, minimum height=2cm] (linear_cell) at (5cm,0cm) {};
        \node[inner sep=0pt, outer sep=0pt, right=0cm of linear_cell] (linear) {};

        \node[cell] (radial_log) at (7.5cm,0.5cm) {$x \mapsto \log |x|$};
        \node[cell] (radial_x) at (7.5cm,-0.5cm) {\sbox0{$x \mapsto \log |x|$}\makebox[\wd0][c]{$x \mapsto x$}};
        \node[inner sep=0pt, outer sep=0pt, right=0cm of radial] (radial) {};

        \node[cell] (linear_2_cell_1) at (10cm,0.5cm) {$x_1$};
        \node[cell] (linear_2_cell_2) at (10cm,-0.5cm) {$x_2$};
        \node[inner sep=0pt, outer sep=0pt, right=0cm of linear_2_cell_1] (linear_2) {};

        \node[cell] (x1_poly_2) at (12.5cm,2cm) {$x_1$};
        \node[cell] (x2_poly_2) at (12.5cm,1cm) {$x_2$};
        \node[cell] (x1x2_poly_2) at (12.5cm,-0cm) {$x_1x_2$};
        \node[cell] (x12_poly_2) at (12.5cm,-1cm) {$x_1^2$};
        \node[cell] (x22_poly_2) at (12.5cm,-2cm) {$x_2^2$};

        \node[cell] (output) at (15cm,0cm) {output};

        \draw[->] (x_in) -- (x_poly);

        \draw[->] (x_poly) -- (4.5cm,0.5cm);
        \draw[->] (x_poly) -- (4.5cm,-0.5cm);

        \draw[->] (5.5cm,0.5cm) -- (radial_log.west);
        \draw[->] (5.5cm,-0.5cm) -- (radial_x.west);

        \draw[->] (radial_log.east) -- (linear_2_cell_1);
        \draw[->] (radial_x.east) -- (linear_2_cell_1);
        \draw[->] (radial_log.east) -- (linear_2_cell_2);
        \draw[->] (radial_x.east) -- (linear_2_cell_2);

        \draw[->] (linear_2_cell_1) -- (x1_poly_2);
        \draw[->] (linear_2_cell_1) -- (x1x2_poly_2);
        \draw[->] (linear_2_cell_2) -- (x1x2_poly_2);
        \draw[->] (linear_2_cell_2) -- (x2_poly_2);
        \draw[->] (linear_2_cell_2) -- (x22_poly_2);
        \draw[->] (linear_2_cell_1) -- (x12_poly_2);

        \draw[->] (x1_poly_2) -- (output);
        \draw[->] (x1x2_poly_2) -- (output);
        \draw[->] (x2_poly_2) -- (output);
        \draw[->] (x22_poly_2) -- (output);
        \draw[->] (x12_poly_2) -- (output);

        \node[inner sep=0pt] (input_label)at (0,-3cm) {Input};
        \node[inner sep=0pt] (polynomial_label) at (2.5cm,-3cm) {Polynomial};
        \node[inner sep=0pt] (linear_label) at (5cm,-3cm) {Linear};
        \node[inner sep=0pt] (radial_label) at (7.5cm,-3cm) {Radial};
        \node[inner sep=0pt] (linear_2_label) at (10cm,-3cm) {Linear};
        \node[inner sep=0pt] (polynomial_2_label) at (12.5cm,-3cm) {Polynomial};
        \node[inner sep=0pt] (output_label) at (15cm,-3cm) {Output};

    \end{tikzpicture}
    \caption{Structure of the PRP model used to learn the Gompertz model equation, in case 2.}
    \label{fig:prp_model_gompertz}
\end{figure}

To begin, we attempt to learn the function using the smallest possible PRP block that can still recover the exact formula (see \cref{fig:prp_model_gompertz}). Nevertheless, this block can generate a large number of functions.

As for the previous application case, the dataset consists of $K=100$ trajectories with $N=500$ constant time steps for each trajectory. Trajectories are observed over the time interval $[0, 2]$. The initial conditions of each trajectory are uniformly sampled between $0$ and $3$.
For the training, we use the same hyperparameters as in the previous case.
The only difference concerns weight initialization, where all the weights randomly initialized according to a normal distribution centered at $0.15$ with a standard deviation of $0.05$. This introduces randomness when running multiple optimizations.

\begin{figure}[htbp]
    \centering
    \begin{subfigure}{\textwidth}
        \includegraphics[width=0.32\linewidth]{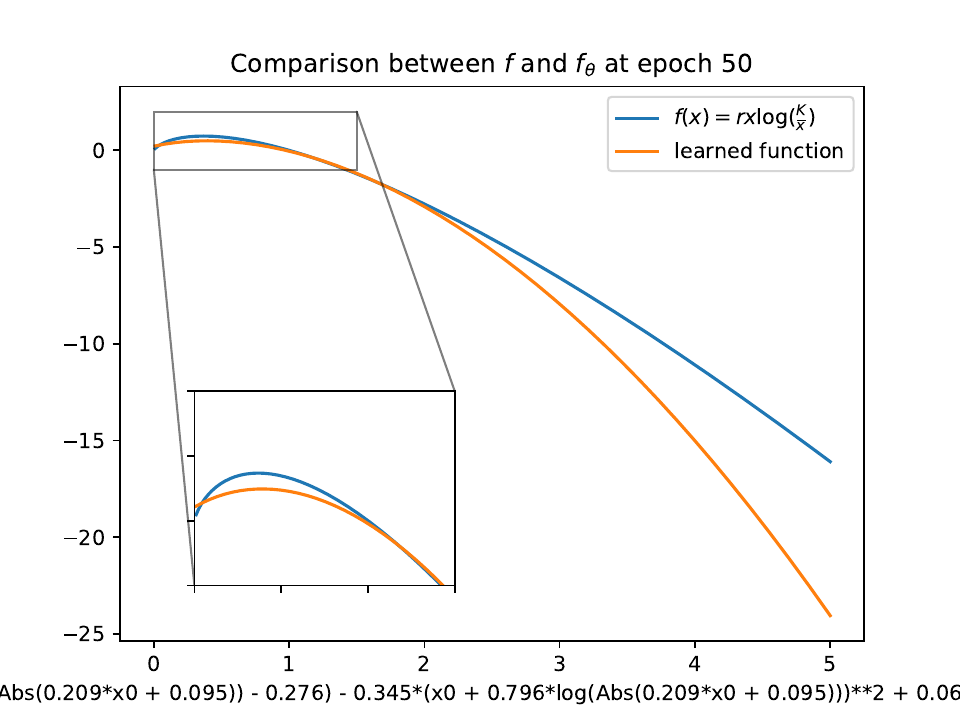}
        \hfill
        \includegraphics[width=0.32\linewidth]{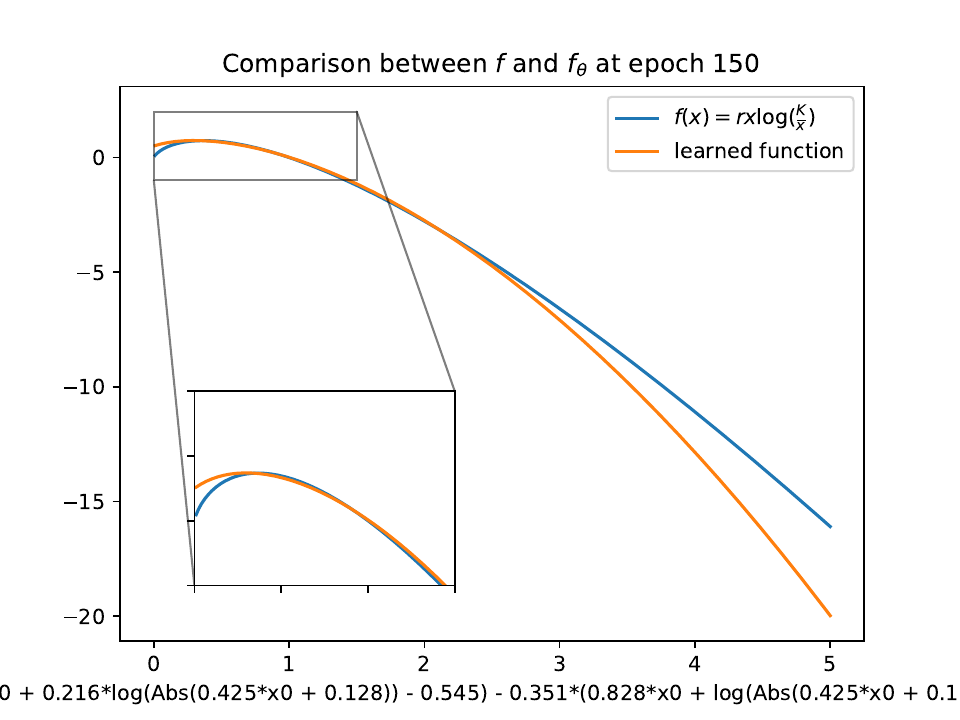}
        \hfill
        \includegraphics[width=0.32\linewidth]{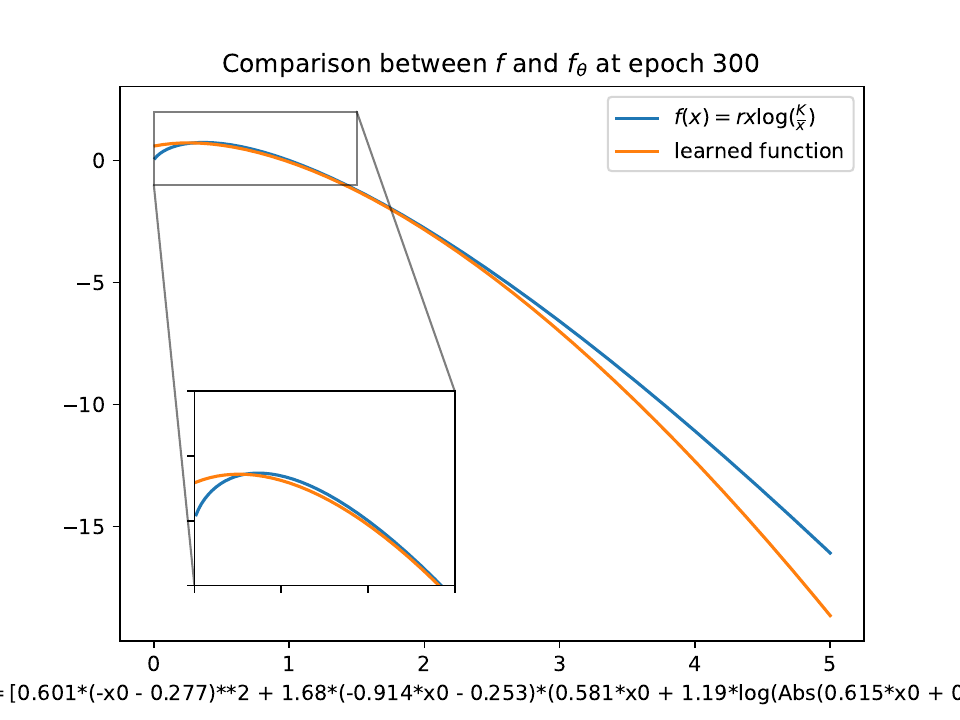}
        \medskip
        \includegraphics[width=0.32\linewidth]{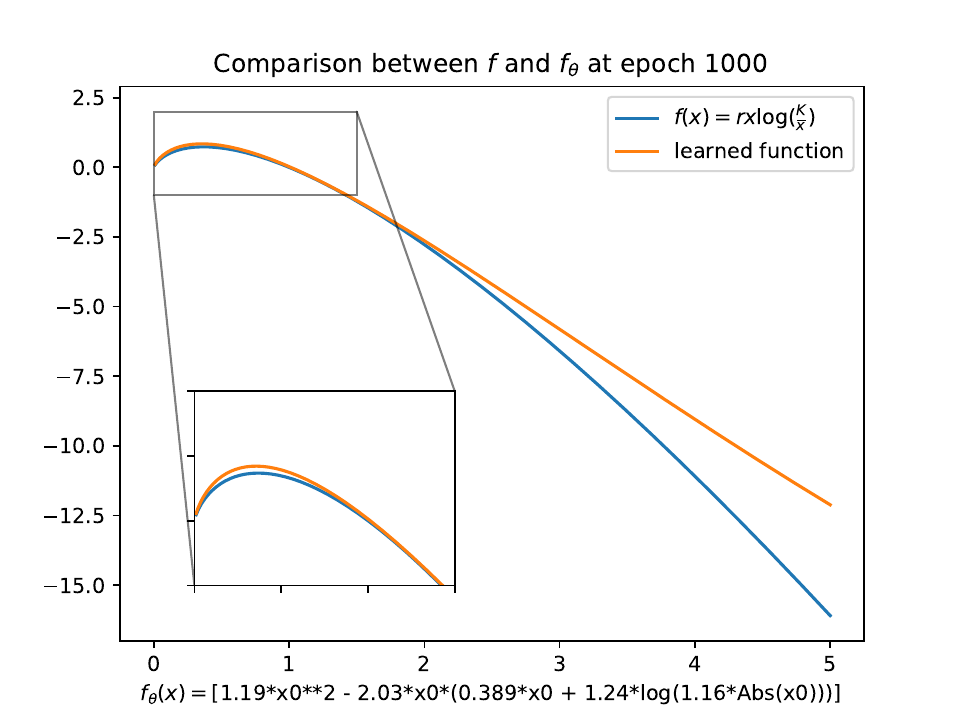}
        \hfill
        \includegraphics[width=0.32\linewidth]{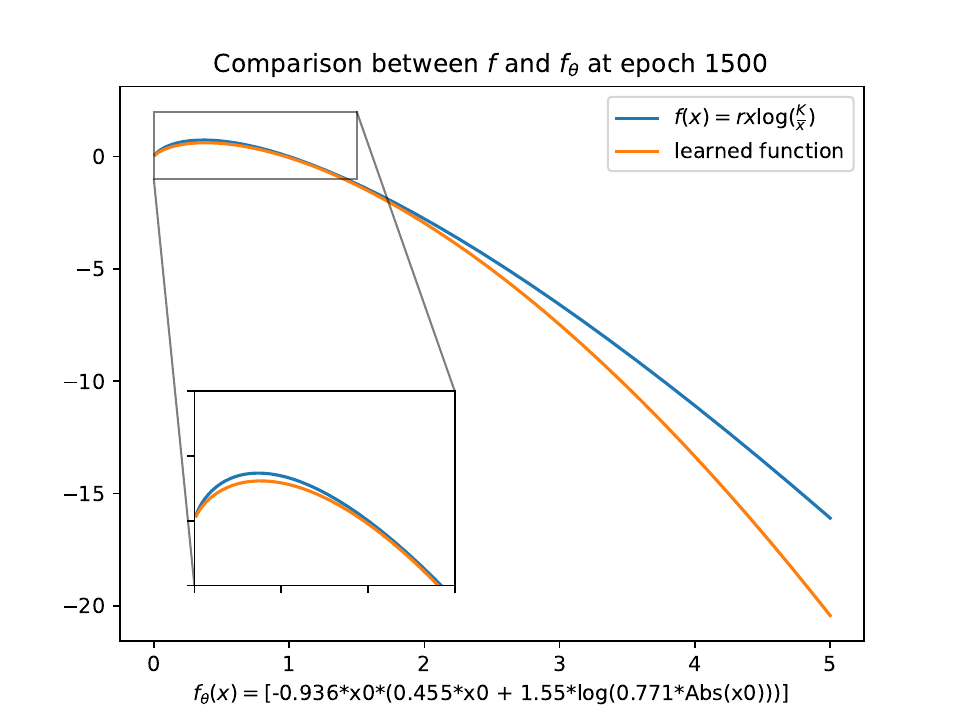}
        \hfill
        \includegraphics[width=0.32\linewidth]{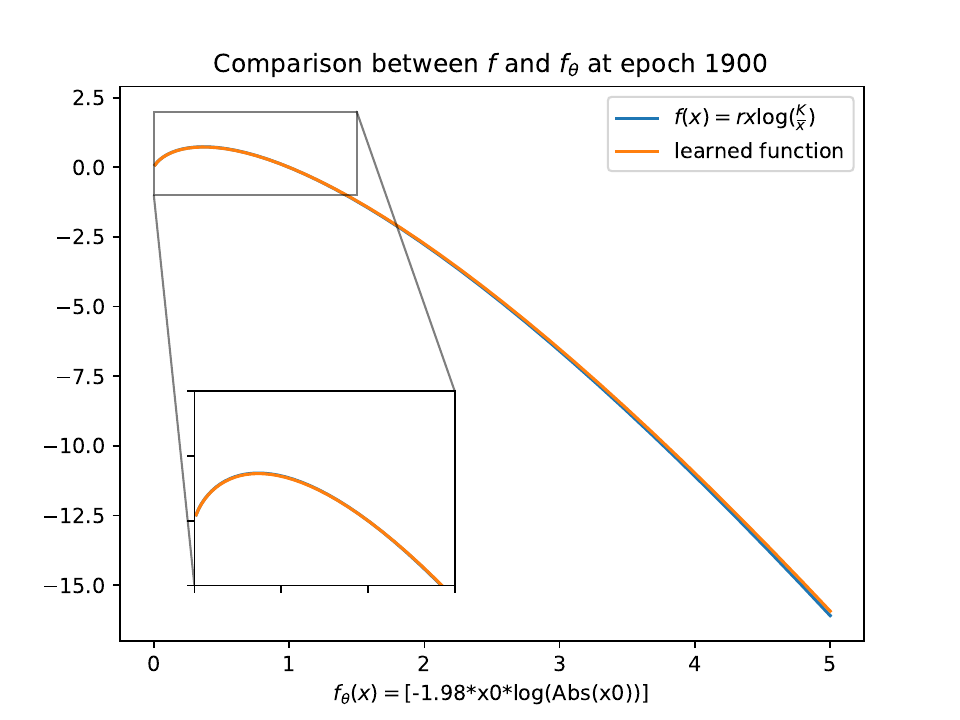}
        \caption{We compare $f$ given by \eqref{eq: Gompertz_speed} (blue curve) and the learned function $f_\theta$ (orange curve), at epochs 50, 150, 300, 1000,  1500, 1900 (from left to right and top to bottom)}
    \end{subfigure}
    \begin{subfigure}{\textwidth}
        \includegraphics[width=0.32\linewidth]{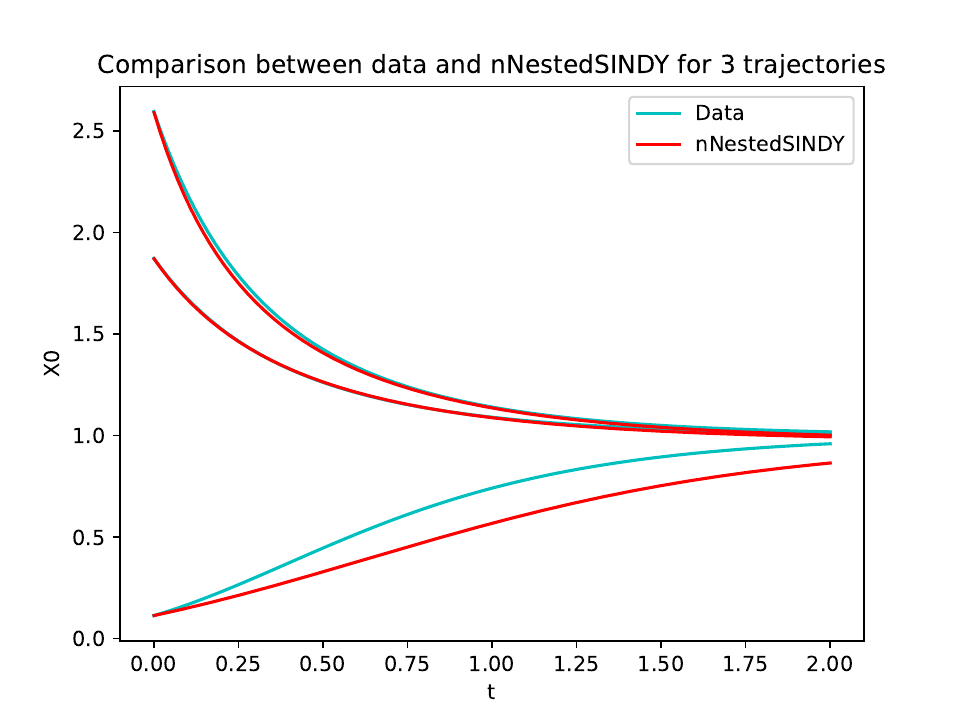}
        \hfill
        \includegraphics[width=0.32\linewidth]{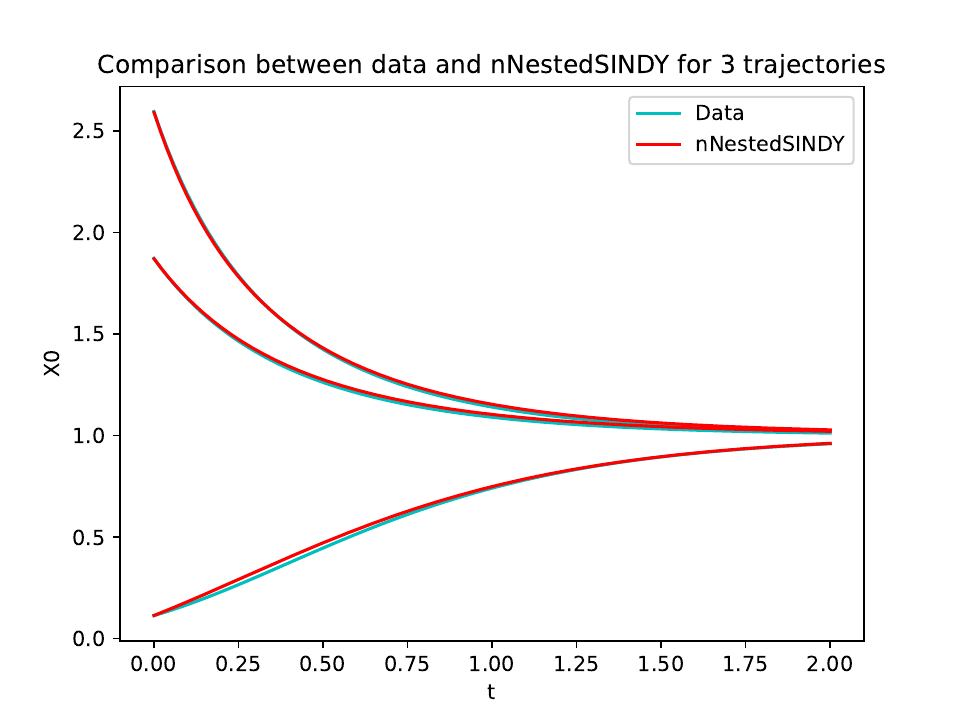}
        \hfill
        \includegraphics[width=0.32\linewidth]{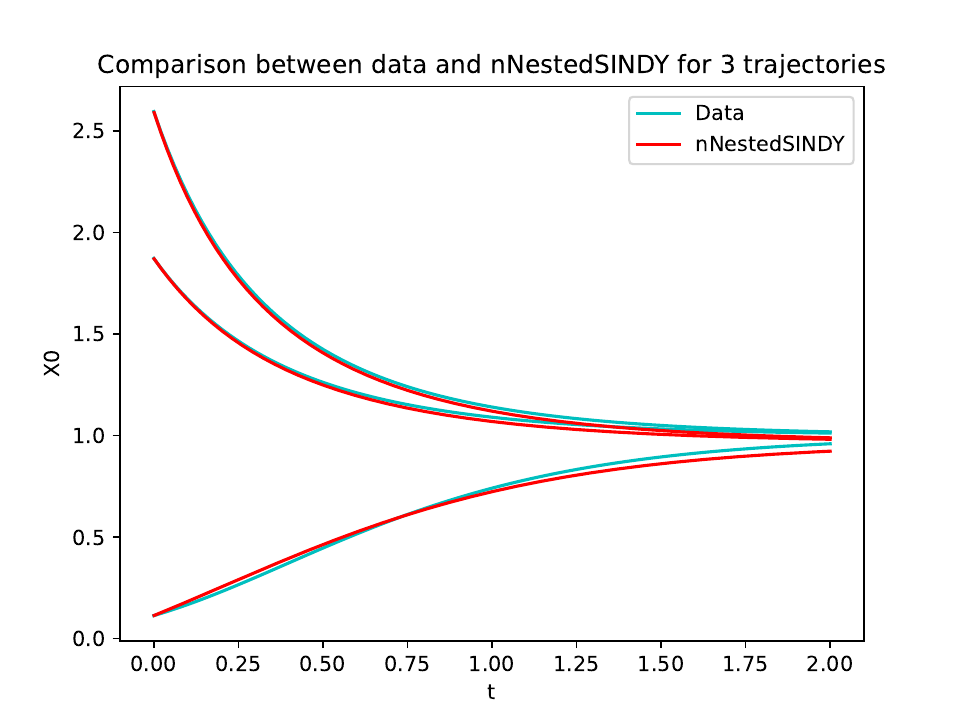}
        \medskip
        \includegraphics[width=0.32\linewidth]{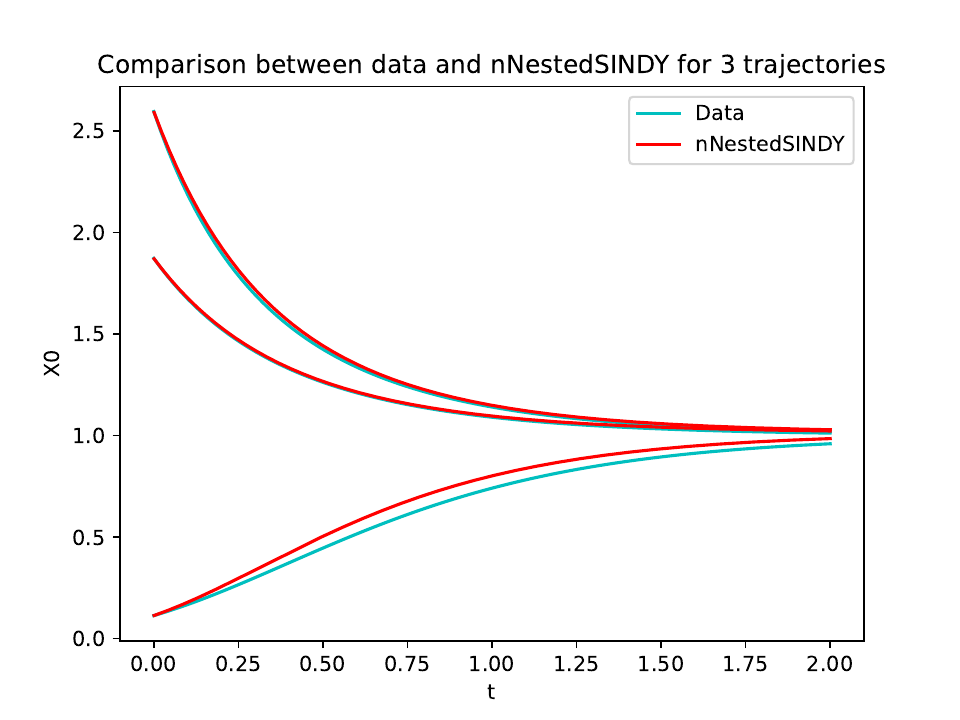}
        \hfill
        \includegraphics[width=0.32\linewidth]{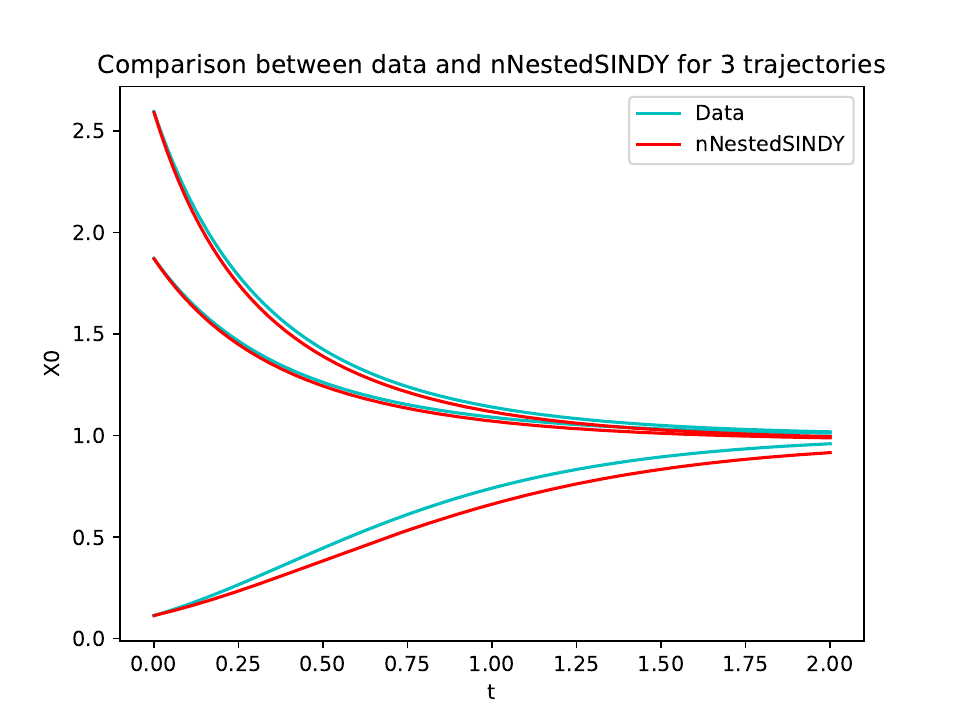}
        \hfill
        \includegraphics[width=0.32\linewidth]{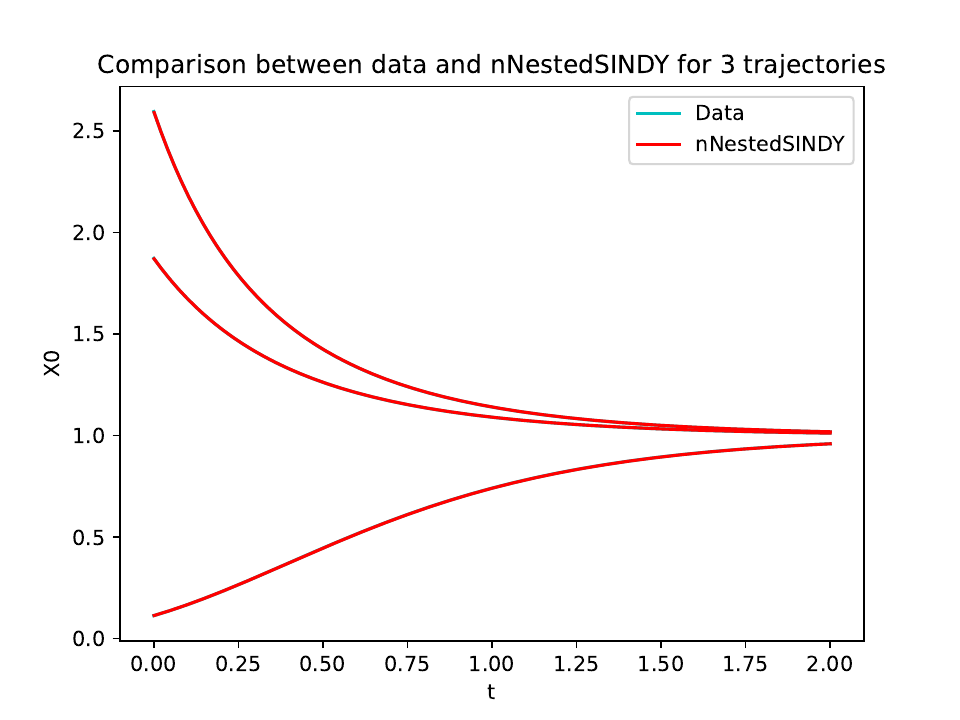}
        \caption{ We compare trajectories of dynamic systems \eqref{eq: edo avec ftheta} (red curves) and \eqref{eq: GompertzODE} (blue curves) for three different initial conditions, at epochs 50, 150, 300, 1000,  1500, 1900 (from left to right and top to bottom)}
    \end{subfigure}
    \caption{Case 2: training of Neural Nested SINDy to recover equation  \eqref{eq: GompertzODE}. The model converges after 1900 epochs.}
    \label{fig: training ode gompertz}
\end{figure}

We train the model for $2000$ epochs. Results are summarized in \cref{fig: training ode gompertz} and \cref{tab:results ode gompertz}. After $1900$ epochs, the model approximates up to an error of $10^{-2}$ the right formula. During the training, the MSE between $f$
and~$f_\theta$ is not always decreasing, but parameters are (slowly) pruned, and finally the model manages to approximate the right formula.


\begin{table}[htbp]
    \caption {Case 2: evolution of the training of Neural Nested SINDy to recover equation \eqref{eq: GompertzODE}. The table contains the formulas of $f_ \theta$ for epochs greater than $1\,000$, the mean squared error between $f_\theta$ and $f$, and the number of pruned parameters in the PRP block, at epochs 50, 150, 300, 1000, 1500, 1900.}
    \label{tab:results ode gompertz}
    \centering
    \begin{tabular}{ccccc}
        \toprule
        Epoch
         & $f_\theta$
         & MSE on $[0,3]$                                           & MSE on $[0,5]$        & \makecell{pruned \\parameters} \\
        \cmidrule(lr){1-5}
        50
         & expression too long to display
         & $1.50 \times 10^{-1}$                                    & $3.83 \times 10^{-1}$ & 4                \\
        150
         & expression too long to display
         & $5.18 \times 10^{-2}$                                    & $1.78 \times 10^{-1}$ & 5                \\
        300
         & expression too long to display

         & $5.31 \times 10^{-2}$                                    & $1.22 \times 10^{-1}$ & 7                \\
        1000
         & $\begin{aligned}
                    1.19x^2 - 2.03x(0.39x + 1.24 \log(1.16\lvert x \rvert ))
                \end{aligned}$

         & $9.17 \times 10^{-2}$                                    & $1.97 \times 10^{-1}$ & 9                \\
        1500
         & $\begin{aligned}
                    -0.94x(0.46x + 1.55 \log(0.77 \lvert x \rvert))
                \end{aligned}$

         & $2.17 \times 10^{-1}$                                    & $3.83 \times 10^{-1}$ & 10               \\
        1900
         & $\begin{aligned}
                    -1.98x\log(\lvert x \rvert))
                \end{aligned}$

         & $1 \times 10^{-2}$                                       & $1 \times 10^{-2}$    & 11               \\
        \bottomrule
    \end{tabular}
\end{table}

As expected (see \cref{sec:example_prp_2D}), PRP models take more effort to retrieve the expected function. In this case we succeed to approximate the true formula of growth speed, but we considered a compact architecture with few dictionary functions, and $20$ times more epochs than in the case of the PR model used in the previous case.
As a consequence, a major way of improving the PRP models lies in improving the optimization process.



\section{Conclusion}
\label{sec:conclusion}

In this study, we explored the capabilities and limitations of the Nested SINDy approach in discovering symbolic representations of dynamical systems from data. Our investigation covered a spectrum of cases, including the identification of ``simple'' functions, as well as the discovery of ordinary differential equations (ODEs) from trajectory data. The versatility of the Nested SINDy framework was demonstrated through various examples, showcasing its ability to approximate complex functions and discover the ODEs underlying data-generating processes.

The Nested SINDy approach extends the original SINDy methodology by incorporating nested structures and neural network architectures, allowing for the identification of complex symbolic expressions that are not directly accessible to traditional methods. This capability is necessary in cases involving compositions and multiplications of functions, where the Nested SINDy method accurately identifies symbolic representations in simple cases, or finds sparse symbolic representations in more complex cases.

Our results confirm that Nested SINDy can recover accurate symbolic expressions for a range of problems, from simple trigonometric functions to the more complex Gompertz model of tumor growth. However, the complexity of the optimization landscape associated with the Nested SINDy approach highlights the need for careful selection of hyperparameters and initialization strategies to ensure approximation of meaningful solutions.

Future work could focus on several areas to enhance the Nested SINDy framework. First, exploring alternative optimization algorithms specifically designed for the unique challenges of nested symbolic regression could improve the efficiency and reliability of the method. Additionally, incorporating mechanisms for automatic selection of the dictionary of basis functions based on preliminary data analysis might streamline the model development process and improve the model's adaptability to different types of dynamical systems.
Furthermore, one could envision integrating a supervised learning component, as in \cite{ParFam}, in which a model is pre-trained to map the relationship between datasets and \sindy-like sparse patterns encoding analytical expressions. This effectively enables the automatic formulation of high-quality initial solutions that can then be refined on a case-by-case basis.
Finally, extending the Nested SINDy approach to handle partial differential equations (PDEs) could open new avenues for discovering the underlying physics of spatially extended systems from data.

In summary, this study underscores the potential of the Nested SINDy approach as a powerful tool for symbolic regression and dynamical system discovery, while also pointing out avenues for further refinement and expansion of the methodology.



\bibliographystyle{plain}
\bibliography{biblio}

\end{document}